\documentclass[reqno,12pt]{article}
\usepackage{a4wide}
\usepackage{amsmath}
\usepackage{amssymb}
\usepackage{amsthm}
\usepackage{mathrsfs}
\usepackage[utf8]{inputenc}
\usepackage{graphicx}

\numberwithin{equation}{section}

\newtheorem{theo}{Theorem}

\newtheorem{lem}{Lemma}
\newtheorem{defi}{Definition}

\theoremstyle{remark}

\def\al{\alpha}

\def\te{\theta}

\def\({\left(}
\def\){\right)}
\def\[{\left[}
\def\]{\right]}

\def\dd{\textup{d}}

\def\Gh{\widehat{\Gamma}}

\def\Qb{\overline{\mathbb Q}}

\newcommand{\N}{\mathbb{N}}
\newcommand{\Z}{\mathbb{Z}}
\newcommand{\Q}{\mathbb{Q}}
\newcommand{\R}{\mathbb{R}}
\newcommand{\C}{\mathbb{C}}
\newcommand{\K}{\mathbb{K}}
\newcommand{\E}{\mathbf{E}}
\newcommand{\Qbar}{\overline{\mathbb Q}}
\newcommand{\etoile}{^*}

\newcommand{\eps}{\varepsilon}

\newcommand{\calD}{{\cal {D}}} 
\renewcommand{\Re}{\textup{Re}\,} 
\renewcommand{\Im}{\textup{Im}\,}

\newcommand{\calFbar}{\overline{\cal F}}
\newcommand{\horsspec}{{\cal S}}
\newcommand{\nvsing}{\Sigma}
\newcommand{\cercle}{{\cal C}(0,R)}
\newcommand{\droite}{d}
\newcommand{\cheminrhoR}{\Gamma_{\rho,R}}
\newcommand{\cheminrhounR}{\Gamma_{\rho_1,R}}
\newcommand{\cheminrhodeR}{\Gamma_{\rho_2,R}}
\newcommand{\cheminrhoinf}{\Gamma_{\rho}}
\newcommand{\groschemin}{\Gamma_R}
\newcommand{\deltarho}{\Delta_\rho}
\newcommand{\pointrho}{z_\rho}
\newcommand{\pointrhoun}{z_{\rho_1}}
\newcommand{\pointrhode}{z_{\rho_2}}
\newcommand{\pointrhop}{z_{\rho_p}}
\newcommand{\Arho}{A_\rho}
\newcommand{\GG}{{\bf G}}

\newcommand{\fct}{\varphi}
\newcommand{\fctalk}{\fct_{\al,k}}

\newcommand{\fctalzero}{\fct_{\al,0}}

\newcommand{\hada}{\star}

\newcommand{\cce}{\omega}
\newcommand{\peal}{\lfloor \al \rfloor}
\newcommand{\listecoeff}{\Phi}
\newcommand{\petiteliste}{\Psi}
\newcommand{\partiexi}{X}
\newcommand{\capa}{\kappa}

\newcommand{\Frac}{{\rm Frac}\,}

\newcommand{\Qbaretoile}{\Qbar\etoile}
\newcommand{\alti}{\widetilde\alpha}
\newcommand{\calK}{{\mathcal K}}
\newcommand{\ddjk}{{\mathfrak D}_{j,k}}

\newcommand{\equi}{\sim}
\newcommand{\devasy}{\approx}
\newcommand{\phiphi}{\phi_{j,s,k}}

\begin{document}

\title{Arithmetic theory of $E$-operators}

\author{S. Fischler and T. Rivoal}
\date{}

\maketitle

\begin{abstract} 
In  [{\em S\'eries Gevrey de type arithm\'etique $I$.
Th\'eor\`emes de puret\'e et de dualit\'e}, Annals of Math. {\bf 151} (2000),   705--740], Andr\'e has introduced $E$-operators, a class of differential operators intimately related to $E$-functions, and constructed local bases of solutions for these operators. In this paper we investigate the arithmetical nature of connexion constants of $E$-operators at finite distance, and of Stokes constants at infinity. We prove that they involve values at algebraic points of $E$-functions in the former case, and in the latter one, values of  $G$-functions and of derivatives of the Gamma function at rational points in a very precise way. As an application, we define and study 
a class of numbers having certain algebraic approximations defined in terms of $E$-functions. These types of approximations are motivated by the convergents to the number $e$, as well as by recent constructions of approximations to Euler's constant and values of the Gamma function. Our results and methods are completely different from those in our paper [{\em On the values of $G$-functions}, Commentarii Math. Helv., to appear], where we have studied  similar questions for $G$-functions.
\end{abstract}

\section{Introduction}

In a seminal paper \cite{YA1}, Andr\'e has introduced $E$-operators, a class of differential operators intimately related to $E$-functions, and constructed local bases of solutions for these operators. In this paper we investigate the arithmetical nature of connexion constants of $E$-operators, and prove that they involve values at algebraic points of $E$-functions or $G$-functions, and values at rational points of derivatives of the Gamma function. As an application, we will focus on algebraic approximations to such numbers, in connection with Aptekarev's famous construction for Euler's constant $\gamma$.

To begin with, let us recall the following definition.

\begin{defi} \label{def:gfunc}
An $E$-function $E$ is a formal power series $E(z)=\sum_{n=0}^{\infty} \frac{a_n}{n!} z^n$ 
such that the coefficients $a_n$ are algebraic numbers and there exists $C>0$ such that:
\begin{enumerate}
\item[$(i)$] the maximum of the moduli of the conjugates of $a_n$ is $\leq C^{n+1}$ for any $n$.

\item[$(ii)$] there exists a sequence of rational integers $d_n$, with $\vert d_n \vert \leq C^{n+1}$, such that 
$d_na_m$ is an algebraic integer for all~$m\le n$.

\item[$(iii)$] $F(z)$ satisfies a homogeneous linear differential equation with 
coefficients in $\Qb(z)$.
\end{enumerate}
\end{defi}
 A $G$-function is defined similarly, as $\sum_{n=0}^{\infty} a_n z^n$ with the same assumptions $(i)$, $(ii)$, $(iii)$; throughout the paper we fix a complex embedding of $\Qbar$.

 We refer to~\cite{YA1} for an overview of the main properties of $E$ and $G$-functions. For the sake of precision, we 
mention that the class of $E$-functions was first defined by Siegel in a more general way, with bounds of the shape $n!^{\eps}$ for 
any $\eps>0$ and any $n\gg _{\eps} 1$, instead of $C^{n+1}$ for all $n\in \N = \{0,1,2,\ldots\}$.  The functions covered by Definition~\ref{def:gfunc} are called $E^*$-functions 
by Shidlovskii~\cite{shid}, and are the ones used in the recent litterature under the denomination $E$-functions (see~\cite{YA1, beukers3, lagarias}); 
 it is believed that both classes coincide.

Examples of $E$-functions include $e^{\alpha z}$ with $\alpha\in\Qbar$, hypergeometric series $_p F_p$ with rational parameters, and Bessel functions. Very precise transcendence (and even algebraic independence) results are known on values of $E$-functions, such as the Siegel-Shidlovskii theorem~\cite{shid}. Beukers' refinement of this result enables one to deduce the following statement 
(see \S\ref{subsecrappelsE}), whose analogue is false for $G$-functions (see~\cite{beukers2} for interesting non-trivial examples): 
 
 \begin{theo}\label{thA}
 An $E$-function with coefficients in a number field $\mathbb K$ 
takes at an algebraic point $\alpha$ either a transcendental 
value or a value in $\mathbb K(\alpha)$.
\end{theo}

 In this paper we consider the following set $\E$, which is analogous to the ring $\GG$ of values at algebraic points of analytic continuations of $G$-functions studied in \cite{gvalues}; we recall that $\GG$ might be equal to $\mathcal P[1/\pi]$, where $\mathcal P$ is the ring of periods (in the sense of Kontsevich-Zagier \cite{KZ}: see \S 2.2 of \cite{gvalues}).

\begin{defi} \label{defi:1} The set $\mathbf E$ 
is defined as the set of all values taken by any $E$-function at any 
algebraic point. 
\end{defi}

Since $E$-functions are entire and $E(\alpha z)$ is an $E$-function for any $E$-function $E(z)$ and any $\alpha\in\Qbar$, we may restrict to values at $z=1$. Moreover $E$-functions form a ring, so that $\E$ is a subring of $\C$. Its group of units contains $\Qbaretoile$ and $\exp(\Qbar)$ because algebraic numbers, $\exp(z)$ and $\exp(-z)$ are $E$-functions. Other elements of $\E$ include values at algebraic points of Bessel functions, and also of any arithmetic Gevrey series of negative order (see \cite{YA1}, Corollaire 1.3.2), for instance Airy's oscillating integral. It seems unlikely that $\E$ is a field and we don't know if we have a full description of its units.

\medskip

A large part of our results is devoted 
to the arithmetic description of {\em connexion constants} or {\em Stokes constants}. 
Any $E$-function 
$E(z)$ satisfies a differential equation $Ly=0$, where $L$ is an $E$-operator (see \cite{YA1}); it is not necessarily minimal and its only possible singularities are
0 and $\infty$. 
Andr\'e has proved~\cite{YA1} that a basis of solutions 
of $L$ at $z=0$ is of the form 
$ 
(E_1(z), \ldots, E_\mu(z))\cdot z^M
$ 
where $M$ is an upper triangular $\mu \times \mu$ matrix with coefficients in $\mathbb Q$ and 
the $E_j(z)$ are $E$-functions. This implies that any local solution $F(z)$ of $L$ at $z=0$ is of the form
\begin{equation}\label{eq:base0}
F(z)=\sum_{j=1}^\mu \Big(\sum_{s \in S_j}
\sum_{k\in K_j} \phiphi z^{s}\log(z)^k\Big) E_j(z) 
\end{equation}
where $S_j \subset \mathbb Q, K_j \subset \mathbb N$ are finite sets and $\phiphi\in \mathbb C$.
Our purpose is to study the connexion constants of $F(z)$, assuming all coefficients $\phiphi$ to be algebraic (with a special focus on the special case where $F(z)$ itself is an $E$-function).

 Any point $\alpha\in\Qbar \setminus\{0\}$ is a regular point of $L$ and there exists a basis of local 
holomorphic solutions $G_1(z), \ldots, G_{\mu}(z)\in\Qbar[[z-\alpha]]$ such that, around $z=\alpha$, 
\begin{equation}\label{eq:EG2}
F(z)=\omega_1G_1(z)+ \cdots + \omega_\mu G_{\mu}(z)
\end{equation}
for some complex numbers $\omega_1, \ldots, \omega_\mu$, called the connexion constants (at finite distance).

\begin{theo}\label{theo:connexfini}
If all coefficients $\phiphi$ in \eqref{eq:base0} are algebraic then 
the connexion constants $\omega_1,\ldots, \omega_\mu$ in~\eqref{eq:EG2} belong to $\E[\log \alpha]$, and even to $\E$ if $F(z)$ is an $E$-function.
\end{theo}

 The situation is much more complicated around $\infty$, which is in general an irregular singularity of $L$; this part is therefore much more involved than the corresponding one for $G$-functions \cite{gvalues} (since $\infty$ is a regular singularity of $G$-operators, the connexion constants of $G$-functions at any $\zeta\in\C\cup\{\infty\}$ always belong to $\GG$). 
 The local solutions at $\infty$ involve divergent series, which give  
 rise to Stokes phenomenon: the expression of an $E$-function $E(z)$ on a given
basis is valid on certain angular sectors, and the connexion constants may change from one sector to another when crossing 
certain rays called anti-Stokes directions. For this reason, we speak of Stokes constants rather than 
connexion constants. 
More precisely, let $\theta\in\R$ and assume that $\theta$ does 
not belong to some explicit finite set (modulo $2\pi$) which contains the anti-Stokes 
directions. Then we compute explicitly the asymptotic expansion
\begin{equation} \label{eqdevintro}
E(z) \devasy \sum_{\rho\in\Sigma} e^{\rho z} 
 \sum_{\al\in S } \sum_{i\in T } \sum_{n=0}^{\infty} c_{\rho, \al,i,n}z^{-n-\al}\log(1/z)^i
 \end{equation}
as $|z| \to \infty$ in a large sector $\theta-\frac{\pi}2-\eps 
\leq \arg(z) \leq \theta+\frac{\pi}2+\eps $ for some $\eps>0$; in precise terms, 
$E(z)$ can be obtained by 1-summation from this expansion (see \S \ref{subsecasyexp}). Here $\Sigma\subset\C$, 
$S\subset\Q$ and $T\subset\N$ are finite subsets, and the coefficients 
$ c_{\rho, \al,i,n}$ are complex numbers; all of them are constructed explicitly 
in terms of the Laplace transform $g(z)$ of $E(z)$, which is annihilated by a 
$G$-operator. In applying or studying~\eqref{eqdevintro} we shall always assume that the sets $\Sigma$, $S$ and $T$ have the least possible cardinality (so that $\alpha-\alpha'\not\in\Z$ for any distinct $\alpha,\alpha'\in S$) and that for any $\alpha$ there exist $\rho$ and $i$ with $ c_{\rho, \al,i,0}\neq 0$. Then the asymptotic expansion~\eqref{eqdevintro} is uniquely determined by $E(z)$ (see \S \ref{subsecasyexp}).

\medskip

One of our main contributions is the value of $ c_{\rho, \al,i,n}$, 
which is given in terms of derivatives of $1/\Gamma$ at $\al\in\Q$ and connexion 
constants of $g(z)$ at its finite singularities~$\rho$. 
Andr\'e has constructed (\cite{YA1}, Th\'eor\`eme 4.3 $(v)$) a basis $H_1(z),\ldots,H_\mu(z)$ 
of formal solutions at infinity of an $E$-operator that annihilates $E(z)$; these solutions involve Gevrey divergent series of order 1, and 
are of the same form as the right hand side of \eqref{eqdevintro}, with algebraic 
coefficients $ c_{\rho, \al,i,n}$. The asymptotic expansion~\eqref{eqdevintro} of $E(z)$ 
in a large sector bisected by $\theta$ can be written on this basis as
\begin{equation} \label{eqintro17}
\cce_{1,\theta} H_1(z) + \ldots + \cce_{\mu,\theta} H_\mu(z)
 \end{equation}
 with Stokes constants $\cce_{i,\theta}$. To identify these constants, we first introduce another important set.

\begin{defi} 
We define ${\bf S}$ as the $\GG$-module generated by all the values of derivatives of the Gamma function at  rational points. It is also the $\GG[\gamma]$-module generated by all the values of $\Gamma$ at  rational points, and it is a ring.
\end{defi} We show in \S\ref{sec:structureS} why the two modules coincide.
The Rohrlich-Lang conjecture (see \cite{YA2} or \cite{MiWperiodes}) implies that the values $\Gamma(s)$, for $s\in\Q$ with $0 < s\leq 1$, are $\Qbar$-linearly independent. We conjecture that these numbers are in fact also $\GG[\gamma]$-linearly independent, so that ${\bf S}$ is the free $\GG[\gamma]$-module they generate.

\medskip

We then have the following result.

\begin{theo} \label{thintrocce}
Let $E(z)$ be an $E$-function, and $\theta\in\R$ be a direction which does not 
belong to some explicit finite set (modulo $2\pi$). 
Then:
\begin{enumerate}
\item[$(i)$] The Stokes constants $\cce_{i,\theta}$ belong to ${\bf S}$.
\item[$(ii)$] All coefficients $ c_{\rho, \al,i,n}$ in~\eqref{eqdevintro} belong to ${\bf S}$.
\item[$(iii)$] Let $\rho\in \Sigma$, $\alpha\in S$, and $n\geq 0$; denote by $k$ the largest $i\in T$ such that $ c_{\rho, \al,i,n}\neq 0$.
 If $k$ exists then for any $i \in T $ the coefficient $ c_{\rho, \al,i,n}$ is a $\GG$-linear combination of $\Gamma(\alpha)$, $\Gamma'(\alpha)$, \ldots, $\Gamma^{(k-i)}(\alpha)$. In particular, $ c_{\rho, \al,k,n} \in \Gamma(\alpha) \cdot \GG$. Here $\Gamma^{(\ell)}(\alpha)$ is understood as $\Gamma^{(\ell)}(1)$ if $\alpha\in \Z_{\leq 0}$. 
 \item[$(iv)$] Let $F(z)$ be a local solution at $z=0$ of an $E$-operator, with algebraic coefficients $\phiphi$ in \eqref{eq:base0}. Then assertions $(i)$ and $(ii)$ hold with $F(z)$ instead of $E(z)$. 
\end{enumerate}
\end{theo}
Assertions $(i)$ and $(iv)$ of Theorem \ref{thintrocce} precise 
Andr\'e's remark in~\cite[p. 722]{YA1}: ``{\em Nous privil\'egierons une approche formelle, qui 
permettrait de travailler sur $\Qbar(\Gamma^{(k)}(a))_{k\in \mathbb N,a\in \mathbb Q}$ plut\^ot que sur 
$\mathbb C$ si l'on voulait}''.

 An important feature of Theorem \ref{thintrocce} (assertion $(iii)$) is that $\Gamma^{(k)}(\alpha)$, for $k\geq 1$ and $\alpha\in\Q\setminus\Z_{\leq 0}$, never appears in the coefficient of a leading term of \eqref{eqdevintro}, but only combined with higher powers of $\log (1/z)$. This motivates the logarithmic factor in \eqref{eq:eapproxgen} below, and explains an observation we had made on Euler's constant: it always appears through $\gamma - \log(1/z)$ (see Eq. \eqref{eqgammalog} in \S \ref{subsecnotationsinf}). Moreover, in $(iii)$, it follows from the remarks made in \S\ref{sec:structureS} that, alternatively, 
$c_{\rho, \al,i,n}= \Gamma(\alpha) \cdot P_{\rho,\al, i,n}(\gamma)$ for some polynomial $P_{\rho,\al, i,n}(X)\in \GG[X]$ of degree $\le k-i$.
 
 The proof of Theorem \ref{thintrocce} is based on Laplace transform, Andr\'e-Chudnovski-Katz's theorem on solutions of $G$-operators, and a specific complex integral (see \cite{YA1}, p. 735). 

\bigskip

As an application of Theorems~\ref{theo:connexfini} and~\ref{thintrocce}, we study sequences of algebraic (or rational) approximations of special interest related to $E$-functions. In \cite{gvalues} we have proved that a complex number $\alpha$ belongs to the fraction field $\Frac \GG$ of $\GG$ if, and only if, there exist sequences $(P_n)$ and $(Q_n)$ of algebraic numbers such that $\lim_n P_n/Q_n =\alpha$ and $\sum_{n\ge 0} P_n z^n$, $\sum_{n\ge 0} Q_n z^n$ are $G$-functions.  We have introduced  this notion in order to give a general framework for irrationality proofs of values of $G$-functions such as zeta values. Such sequences are called $G$-approximations of $\alpha$, when 
$P_n$ and $Q_n$ are rational numbers. We drop this last assumption in the context of $E$-functions (see \S \ref{subsecrappelsE}), and consider the following definition.
 
\begin{defi}\label{def:Eapprox}
Sequences $(P_n)$ and $(Q_n)$ of algebraic numbers are said to be {\em $E$-approxima\-tions} of $\alpha\in\C$ if 
$$
\lim_{n\to +\infty} \frac{P_n}{Q_n} =\alpha
$$
and 
$$
\sum_{n=0}^{\infty} P_n z^n= A(z)\cdot E\big(B(z)\big), \quad \sum_{n=0}^{\infty} Q_n z^n= C(z)\cdot F\big(D(z)\big)
$$
where $E$ and $F$ are $E$-functions, $A, B, C, D$ are algebraic functions in $\Qbar[[z]]$ 
with $B(0)=D(0)=0$.
\end{defi}
This definition is motivated by the fact that many sequences of approximations to classical numbers are $E$-approximations, for instance diagonal Pad\'e approximants to $e^z$ and in particular the convergents of the continued fraction expansion of $e$ (see \S \ref{ssec:example}). Elements in $\Frac \GG$ also have $E$-approximations, since $G$-approximations $(P_n)$ and $(Q_n)$ of a complex number always provide $E$-approximations $P_n/n!$ and $Q_n/n!$ of the same number. In \S \ref{ssec:example}, we construct  $E$-approximations to $\Gamma(\alpha)$ for any $\alpha\in\Q\setminus\Z_{\leq 0}$, $\alpha<1$, by letting 
$
E_\alpha(z)=\sum_{n=0}^{\infty} \frac{z^n}{n!(n+\alpha)}
$, $Q_n(\alpha)=1$, and defining 
$P_n(\alpha)$ by the series expansion (for $\vert z\vert <1$)
$$
\frac1{(1-z)^{\alpha+1}} E_\alpha\left(-\frac{z}{1-z}\right) =\sum_{n=0}^{\infty} P_n(\alpha) z^n \in \mathbb Q[[z]];
$$
then $\lim_n P_n(\alpha) = \Gamma(\alpha)$. The number $\Gamma(\alpha)$ appears in this setting as a Stokes constant. The condition $\alpha<1$ is harmless because we 
readily deduce  $E$-approximations to $\Gamma(\alpha)$ for any $\alpha\in\Q$, $\alpha>1$, by means of the functional equation $\Gamma(s+1)=s\Gamma(s)$. Moreover, since $\frac1{(1-z)^{\alpha+1}} E_\alpha\left(-\frac{z}{1-z}\right)$ is holonomic, the sequence $(P_n(\alpha))$ satisfies a linear recurrence, of order $3$ with polynomial coefficients in $\mathbb Z[n,\alpha]$ of total degree $2$ in $n$ and $\alpha$; see \S\ref{ssec:example}. This construction is simpler than that in~\cite{rivoal3} but the convergence to $\Gamma(\alpha)$ is slower.

\medskip

Definition \ref{def:Eapprox} enables us to consider an interesting class of numbers: those having $E$-approximations. Of course this is a countable subset of $\C$. We have seen that it contains all values of the Gamma function at rational points $s$, which are conjectured to be irrational if $s\not\in\Z$; very few results are known in this direction (see~\cite{MiWperiodes}), and using suitable $E$-approximations may lead to prove new ones. 

However we conjecture that Euler's constant $\gamma$ does not have $E$-approximations: all approximations we have thought of seem to have generating functions not as in Definition~\ref{def:Eapprox}. 
This is a reasonable conjecture in view of Theorem \ref{theo:eapprox} we are going to state now.

 Given two subsets $X$ and $Y$ of $\mathbb C$, we set 
$$
X\cdot Y=\big\{xy\,\big\vert\, x\in X, y\in Y\big\}, \quad 
\displaystyle \frac{X}{Y}=\Big\{\frac xy \,\Big\vert \,x\in X, y\in Y\setminus\{0\}\Big\}.
$$ 
We also set 
$\Gamma({\mathbb Q})=\{\Gamma(x)\vert x\in \mathbb Q\setminus \mathbb Z_{\le 0}\}$. 
If $X$ is a ring then we denote by $\Frac X = \frac{X}{X}$ its field of fractions. We recall \cite{gvalues} that $B(x,y)$ belongs to the group of units $\GG \etoile$ of $\GG$ for any $x,y\in\Q$, so that $\Gamma$ induces a group homomorphism $\Q \to \C\etoile/\GG\etoile$ (by letting $\Gamma(-x) = 1$ for $x\in\N$). Therefore $\Gamma(\Q) \cdot \GG\etoile$ is a subgroup of $\C\etoile$, and so is $\Gamma(\Q) \cdot \exp(\Qbar) \cdot \Frac\GG$; for future reference we write
\begin{equation} \label{eq123}
 \Gamma(\Q) \cdot\Gamma(\Q) \subset \Gamma(\Q)\cdot \GG \quad \mbox{ and } \quad \frac{\Gamma(\Q) }{\Gamma(\Q) } \subset \Gamma(\Q) \cdot \GG.
 \end{equation}

\begin{theo}\label{theo:eapprox} The set of numbers having $E$-approximations contains 
\begin{equation}\label{eq:subset1}
\frac{ \E \cup \Gamma(\Q)}{ \E \cup \Gamma(\Q)} \cup \Frac\GG
\end{equation}
and it is contained in 
\begin{equation}\label{eq:subset2}
\frac{ \E \cup (\Gamma(\Q) \cdot \GG)}{ \E \cup (\Gamma(\Q) \cdot \GG)} \cup \Big(\Gamma(\Q) \cdot \exp(\Qbar) \cdot \Frac\GG \Big).
\end{equation}
\end{theo}
The proof of \eqref{eq:subset1} is constructive; the one of \eqref{eq:subset2} is based on an explicit determination of the asymptotically dominating term of a sequence $(P_n)$ as in Definition \ref{def:Eapprox}. This determination is based on analysis of singularities, the saddle point method, asymptotic expansions \eqref{eqdevintro} of $E(z)$, and Theorems \ref{theo:connexfini} and \ref{thintrocce}; it is of independent interest (see 
Theorem \ref{theoasypn} in \S \ref{sec:asympPn}). The dominating term comes from the local behaviour of $E(z)$ at some $z_0\in\C$ (providing elements of $\E$, in connection with Theorem \ref{theo:connexfini}) or at infinity (providing elements of $\Gamma(\Q)\cdot\GG$; Theorem \ref{thintrocce} is used in this case). This dichotomy leads to the unions in \eqref{eq:subset1} and \eqref{eq:subset2}; it makes it unlikely for 
 the set of numbers having $E$-approximations to be a field, or even a ring. We could have obtained a field by restricting Definition \ref{def:Eapprox} to the case where $B(z) = D(z) = z$ and $A(z)$, $C(z)$ are not polynomials, since in this case the behavior of $E(z)$ at $\infty$ would not come into the play; this field would be simply $\Frac\E$. 

\medskip
 
 It seems likely that there exist numbers having $E$-approximations but no $G$-approxi\-ma\-tions, because 
conjecturally $\Frac\E \cap \Frac\GG=\Qbar$ and $\Gamma(\Q)\cap\Frac\GG = \Q$. It is also an open question to prove that the number $\Gamma^{(n)}(s)$ does not have $E$-approximations, for $n\geq 1$ and $s\in \Q\setminus \mathbb Z_{\le 0}$. 
To obtain approximations to these numbers, one can consider the following generalization of Definition \ref{def:Eapprox}: we replace $A(z)\cdot E(B(z))$ (and also $C(z)\cdot F(D(z))$) with a finite sum
\begin{equation}\label{eq:eapproxgen}
\sum_{i,j,k,\ell} \alpha_{i,j,k,\ell} \log(1-A_{i}(z))^j \cdot B_{k}(z) \cdot E_\ell\big(C(z)\big)
\end{equation}
where $\alpha_{i,j,k,\ell}\in \Qbar$, 
$A_i(z), B_k(z), C (z)$ are algebraic functions in $\Qbar[[z]]$, $A_i(0)=C (0)=0$, and $E_\ell(z)$ are $E$-functions. 
For instance, let us consider the $E$-function
$
E(z)=\sum_{n=1}^{\infty} \frac{z^n}{n!n}
$
and define $P_n$ by the series expansion
 (for $\vert z\vert <1$)
\begin{equation} \label{eq44}
\frac{\log(1-z)}{1-z}-\frac{1}{1-z} E\Big(-\frac{z}{1-z}\Big) = \sum_{n=0}^{\infty} P_n z^n\in \mathbb Q[[z]].
\end{equation}
Then we prove in \S \ref{ssec:extended} that 
$\lim_n P_n=\gamma$, so that letting $Q_n=1$ we obtain $E$-approximations of Euler's constant in this extended sense. Since $\frac{\log(1-z)}{1-z}-\frac1{1-z} E\left(-\frac{z}{1-z}\right)$ is holonomic, the sequence $(P_n)$ satisfies a linear recurrence, of order $3$ with polynomial coefficients in $\mathbb Z[n]$ of degree $2$; see \S\ref{ssec:extended}. Again, this construction is much simpler than those in~\cite{Aptekarev, HP2, rivoal1} but the convergence to $\gamma$ is slower. A construction similar to~\eqref{eq44}, based on an immediate generalization of the final equation for $\Gamma^{(n)}(1)$ in~\cite{Michigan}, shows that the numbers $\Gamma^{(n)}(s)$ have $E$-approximations in the extended sense of~\eqref{eq:eapproxgen} for any integer $n\ge 0$ and any rational number $s\in \mathbb Q\setminus \mathbb Z_{\le 0}$.

The set of numbers having such approximations is still countable, and we prove in \S \ref{ssec:extended} that it is contained in 
\begin{equation}\label{eq:subset3}
\frac{( \E \cdot \log(\Qbaretoile)) \cup {\bf S} }{( \E \cdot \log(\Qbaretoile)) \cup {\bf S} } \cup \Big( \exp(\Qbar) \cdot \Frac{\bf S} \Big) 
\end{equation}
where $\log(\Qbaretoile) = \exp^{-1}(\Qbaretoile)$. 

The generalization~\eqref{eq:eapproxgen} does not cover all interesting constructions of approximations to derivatives of Gamma values in the literature. 
For instance, it does not seem that Aptekarev's or the second author's approximations to $\gamma$ (in~\cite{Aptekarev} and~\cite{rivoal1} respectively) can be described by~\eqref{eq:eapproxgen}. This is also not the case of Hessami-Pilehrood's approximations to $\Gamma^{(n)}(1)$ in \cite{HP, HP2} but in certain cases their generating functions involve sums of products of $E$-functions at various algebraic functions, rather linear forms in $E$-functions at one algebraic function as in~\eqref{eq:eapproxgen}. Another possible generalization of~\eqref{eq:eapproxgen} is to let $\alpha_{i,j,k,\ell}\in \E$; 
we describe such an example in \S\ref{ssec:extended}, related to the continued fraction $[0;1,2,3,4,\ldots]$ whose partial quotients are the consecutive positive integers.

\bigskip

The structure of this paper is as follows. In \S\ref{sec:structureS}, we discuss the properties of ${\bf S}$. In \S \ref{sec2} we prove our results at finite distance, namely Theorems \ref{thA} and 
\ref{theo:connexfini}. Then we  discuss in \S \ref{subsecasyexp} the definition and basic properties of asymptotic expansions. This allows us to prove Theorem \ref{thintrocce} in \S \ref{sec3}, and to determine in \S \ref{sec:asympPn} the asymptotic behavior of sequences $(P_n)$ as in Definition~\ref{def:Eapprox}. Finally, we gather in \S \ref{sec:eapprox} all results related to $E$-approximations.

\section{Structure of ${\bf S}$}\label{sec:structureS}

In this short section, we discuss the structural  properties of the $\GG$-module ${\bf S}$ generated by the numbers $\Gamma^{(n)}(s)$, for $n\ge0$, $s\in \mathbb Q\setminus \mathbb Z_{\le 0}$. It is not used in the proof of our theorems. 

The Digamma function $\Psi$ is defined as the logarithmic derivative of the Gamma function. We have
$$
\Psi(x)=-\gamma+\sum_{k=0}^{\infty}\Big(\frac{1}{k+1}-\frac{1}{k+x}\Big) \quad \textup{and} \quad \Psi^{(n)}(x)
=\sum_{k=0}^{\infty} \frac{(-1)^{n+1} n!}{(k+x)^{n+1}}\quad (n\ge 1).
$$
From the relation 
$\Gamma'(x)=\Psi(x)\Gamma(x)$, we can prove by induction on the integer $n\ge 0$ that 
$$
\Gamma^{(n)}(x)= \Gamma(x)\cdot P_n\big(\Psi(x),\Psi^{(1)}(x),\ldots, \Psi^{(n-1)}(x)\big)
$$
where $P_n(X_1, X_2, \ldots, X_n)$ is a polynomial with integer coefficients. Moreover, the term of maximal degree in $X_1$ 
is $X_1^n$.

It is well-known that $\Psi(s)\in -\gamma + \GG$ (Gauss' formula,~\cite[p.~13, Theorem 1.2.7]{Andrews}) and that 
$\Psi^{(n)}(s)\in \GG$ for any $n\ge 1$ and any $s\in\Q\setminus \mathbb{Z}_{\le 0}$. It follows that 
\begin{equation}\label{eq:deriGamma}
\Gamma^{(n)}(s)=\Gamma(s)\cdot P_n\big(\Psi(s),\Psi^{(1)}(s),\ldots, \Psi^{(n-1)}(s)\big)=
\Gamma(s)\cdot Q_{n,s}(\gamma)
\end{equation}
where $Q_{n,s}(X)$ is a polynomial with coefficients in $\GG$, of degree $n$ and leading coefficient equal to $(-1)^n$.

We are now ready to prove that ${\bf S}$ coincides with the $\GG[\gamma]$-module $\widehat{{\bf S}}$ generated by the numbers 
$\Gamma(s)$, for $s\in \mathbb Q\setminus \mathbb Z_{\le 0}$. Indeed, Eq.~\eqref{eq:deriGamma} shows immediately that ${\bf S}\subset \widehat{{\bf S}}$. For the converse inclusion $\widehat{{\bf S}} \subset {\bf S}$, it is enough to show that $\Gamma(s)\gamma^n\in {\bf S}$ for any $n\ge0$, $s\in \mathbb Q\setminus \mathbb Z_{\le 0}$. This can be proved by induction on $n$ from~\eqref{eq:deriGamma} because we can rewrite it as
$$
\Gamma(s)\gamma^n = (-1)^n \Gamma^{(n)}(s)+\Gamma(s)\cdot \widehat{Q}_{n,s}(\gamma)
$$
for some polynomial $\widehat{Q}_{n,s}(X)$ with coefficients in $\GG$ and degree $\le n-1$.

The module $\widehat{{\bf S}}$ is easily proved to be a ring. Indeed, defining Euler's Beta function  $B(x,y)=\frac{\Gamma(x)\Gamma(y)}{\Gamma(x+y)}$, then for any $x,y\in\Q\setminus \mathbb Z_{\le 0}$ we have $\Gamma(x)\Gamma(y)=\Gamma(x+y)B(x,y)\in \widehat{{\bf S}}$ because $B(x,y)\in \GG$ (see \cite{gvalues}). This can also be proved directly 
from the definition of ${\bf S}$: for any $x,y\in\Q\setminus \mathbb Z_{\le 0}$, we have
\begin{eqnarray*}
\Gamma^{(m)}(x) \Gamma^{(n)}(y) 
&=& \frac{\partial^{m+n}}{\partial x^m \partial y^n} \Gamma(x+y) B(x,y) \\
&=& \sum_{i=0}^m \sum_{j=0}^n \binom{m}{i} \binom{n}{j} \Gamma^{(i+j)}(x+y) \frac{\partial^{m+n-i-j}}{\partial x^{m-i}\partial y^{n-j}} B(x,y) \in {\bf S}
\end{eqnarray*}
because  $ \frac{\partial^{m+n-i-j}}{\partial x^{m-i}\partial y^{n-j}} B(x,y) \in \GG$, arguing as in \cite{gvalues} for the special case $m-i = n-j = 0$.

\section{First results on values of $E$-functions} \label{sec2}
 
\subsection{Around Siegel-Shidlovskii and Beukers' theorems} \label{subsecrappelsE}

To begin with, let us mention the following result. It is proved in~\cite{gvalues} (and due to the referee of that paper) in the case $\mathbb K=\mathbb Q(i)$; actually the same proof, which relies on Beukers' version~\cite{beukers3} of the Siegel-Shidlvoskii theorem, works for any $\mathbb K$. 

\begin{theo} \label{propefcts}
Let $E(z)$ be an $E$-function with coefficients in some number field $\mathbb K$, and $\al,\beta \in \Qbar$ be 
such that $E(\al) = \beta $ or $E(\al) = e^\beta $. Then $\beta\in\mathbb K(\al)$.
\end{theo}

This result implies Theorem \ref{thA} stated in the introduction;  without further hypotheses $E(\alpha)$ may really belong to $\K(\alpha)$, since  if $E(z)$ is an $E$-function then so is $(z-\alpha)E(z)$.

\bigskip

Theorem \ref{propefcts} shows that if we restrict the coefficients of $E$-functions to a given number field then the set of values we obtain is a proper subset of $\E$. In this respect the situation is completely different from the one with $G$-functions, since any element of $\GG$ can be written \cite{gvalues} as $f(1)$ for some $G$-function $f$ with Taylor coefficients in $\Q(i)$. This is also the reason why we did not restrict to rational numbers $P_n$, $Q_n$ in Definition \ref{def:Eapprox}.

\subsection{Connexion constants at finite distance}\label{ssec:11}

Let us prove Theorem \ref{theo:connexfini} stated in the introduction; the strategy 
 is analogous to the corresponding one with $G$-functions \cite{gvalues}, and even easier because $E$-functions are entire.

\begin{proof}
We write 
$$
L= \frac{d^\mu}{d z^\mu} + a_{\mu-1}(z) \frac{d^{\mu-1}}{d z^{\mu-1}}
+ \cdots +a_1(z) \frac{d}{dz} + a_0(z),
$$
where $a_j\in \Qbar (z)$. Then $z=0$ is the 
only singularity at finite distance of $L$, and it is a regular singularity with rational 
exponents (see \cite{YA1}). Hence, any wronskian $W(z)$ of $L$, i.e. any solution of the differential 
equation $y'(z)+a_{\mu-1}(z)y(z)=0$, is of the form
$ 
W(z)= c z^\rho
$ 
with $c\in \mathbb C$ and $\rho \in \mathbb Q$. We denote by $W_G(z)$ the wronskian 
of $L$ built on the functions $G_1(z), \ldots, G_\mu(z)$: 
$$
W_G(z)= \left\vert 
\begin{matrix}
 G_{1}(z) &\cdots &G_{\mu}(z)
\\
G_{1}^{(1)}(z) &\cdots 
&G_{\mu}^{(1)}(z)
\\
\vdots &\cdots &\vdots
\\
G_{1}^{(\mu-1)}(z) &\cdots &G_{\mu}^{(\mu-1)}(z)
\end{matrix}
\right\vert.
$$
All  functions $G_j ^{(k)}(z)$ are holomorphic at $z=\alpha$ with Taylor coefficients 
in $\Qbar$. Hence, $W_G(\alpha) \in \Qbar$ and is non zero because we also have 
$W_G(\alpha)=c \alpha^\rho$ for some $c$, with $c\neq 0$ because the $G_j$ form a basis of solutions of $L$.

We now differentiate \eqref{eq:EG2} to obtain the relations 
$$F^{(k)}(z)=\sum_{j=1}^{\mu} \omega_j G_j^{(k)}(z), \quad k=0,\ldots, \mu-1$$
for any $z$ in some open disk $\mathcal D$ centered at $z=\alpha$. We interpret these equations as a linear system with unknowns  $\omega_j$, and   solve  it using Cramer's rule. We obtain this way that
\begin{equation}\label{eq:connectiondeter}
\omega_j=\frac1{W_G(z)} 
 \left\vert 
\begin{matrix}
 G_{1}(z) &\cdots &G_{j-1}(z)&F(z) 
&G_{j+1}(z)&\cdots &G_{\mu}(z)
\\
G_{1}^{(1)}(z) &\cdots &G_{j-1}^{(1)}(z)
&F^{(1)}(z) &G_{j+1}^{(1)}(z)&\cdots 
&G_{\mu}^{(1)}(z)
\\
\vdots &\cdots & \vdots & \vdots & \vdots&\cdots&\vdots
\\
G_{1}^{(\mu-1)}(z) &\cdots &G_{j-1}^{(\mu-1)}(z)
&F^{(\mu-1)}(z) &G_{j+1}^{(\mu-1)}(z)
&\cdots &G_{\mu}^{(\mu-1)}(z)
\end{matrix}
\right\vert 
\end{equation}
for any $z \in \mathcal D\setminus\{0\}$, since $W_G(z)\neq 0$.

We now choose $z=\alpha$. As already said, $1/W_G(\alpha), G_j^{(k)}(\alpha)\in \Qbar \subset \E$. If 
we assume that $F(z)$ is an $E$-function, this is also the case of its derivatives, so 
that $F^{(k)}(\alpha) \in \E$ for all $k\ge 0$ and \eqref{eq:connectiondeter} 
implies that $\omega_j \in \E$.
To prove the general case, we simply observe that if $F(z)$ is given by \eqref{eq:base0} with algebraic coefficients $\phiphi$  then all 
 derivatives of $F(z)$ 
at $z=\alpha$ belong to $ \E[\log(\alpha)]$.
\end{proof}

\section{Stokes constants of $E$-functions} \label{sec3}

In this section we construct explicitly the asymptotic expansion of an $E$-function: our main result is Theorem \ref{theoprecis}, stated in \S \ref{subsecnotationsinf} and proved in \S \ref{subsecdemtheoprecis}. Before that we discuss in \S \ref{subsecasyexp} the asymptotic expansions used in this paper. Finally we show in \S \ref{subsec34} that Theorem \ref{theoprecis} implies Theorem \ref{thintrocce}.

Throughout this section, we let $\Gh:=1/\Gamma$ for simplicity.

\subsection{Asymptotic expansions} \label{subsecasyexp}

The asymptotic expansions used throughout this paper are defined as follows.

\begin{defi} \label{defiasy}
Let $\theta\in\R$, and $\Sigma\subset\C$, $S\subset\Q$, $T\subset\N$ be finite subsets. Given
complex numbers $ c_{\rho, \al,i,n}$, we write
\begin{equation} \label{eqasy1}
f(x) \devasy \sum_{\rho\in\Sigma} e^{\rho x} 
 \sum_{\al\in S } \sum_{i\in T } \sum_{n=0}^{\infty} c_{\rho, \al,i,n}x^{-n-\al}(\log(1/x))^i
 \end{equation}
and say that the right hand side is the asymptotic expansion of $f(x)$ in a large sector bisected by 
the direction $\theta$, if there exist $\eps, R, B, C > 0$ and, for any $\rho\in\Sigma$, a function $f_\rho(x)$ holomorphic on
$$U = \Big\{x\in\C, \,\,|x|\geq R, \, \, \theta-\frac{\pi}2-\eps \leq \arg(x) \leq \theta+\frac{\pi}2+\eps \Big\},
$$
such that 
$$f(x) = \sum_{\rho\in\Sigma} e^{\rho x} f_\rho(x)
 $$
 and
$$
\Big| f_\rho(x) - \sum_{\al\in S } \sum_{i\in T } \sum_{n=0}^{N-1} c_{\rho, \al,i,n}x^{-n-\al}(\log(1/x))^i\Big| \leq C^N N! |x|^{B-N} 
$$
for any $x\in U$ and any $N\geq 1$.
\end{defi}

This means exactly (see \cite[\S\S 2.1 and 2.3]{Ramis})  that for any $\rho\in\Sigma$, 
\begin{equation} \label{eqasy2}
 \sum_{\al\in S } \sum_{i\in T } \sum_{n=0}^{N-1} c_{\rho, \al,i,n}x^{-n-\al}(\log(1/x))^i
 \end{equation}
 is 1-summable in the direction $\theta$ and its sum is $f_\rho(x)$. In particular, using a result of
 Watson (see \cite[\S 2.3]{Ramis}), the sum $f_\rho(x)$ is determined by its expansion \eqref{eqasy2}. 
Therefore the asymptotic expansion on the right hand side of \eqref{eqasy1} determines the function $f(x)$ 
(up to analytic continuation). The converse is also true, as the following lemma shows.
 
 \begin{lem} \label{lemasyunique}
 A given function $f(x)$ can have at most one asymptotic expansion in the sense of Definition \ref{defiasy}.
\end{lem}

Of course we assume implicitly in Lemma \ref{lemasyunique} (and very often in this paper) that $\Sigma$, $S$ and $T$ in \eqref{eqasy1} cannot trivially be made smaller, and that for any $\alpha$ there exist $\rho$ and $i$ with $c_{\rho, \al,i,0}\neq 0$. 

\begin{proof}
We proceed by induction on the cardinality of $\Sigma$. If the result holds for proper subsets of 
$\Sigma$, we choose $\theta'$ very close to $\theta$ such that the complex numbers $\rho e^{i\theta'}$, 
$\rho\in\Sigma$, have pairwise distinct real parts and we denote by $\rho_0$ the element of $\Sigma$ for
 which $\Re( \rho_0 e^{i\theta'})$ is maximal. Then the asymptotic expansion \eqref{eqasy2} of $f_{\rho_0}(x)$ 
is also an asymptotic expansion of $e^{-\rho_0 x} f(x)$ as $|x|\to\infty$ with $\arg(x) = \theta'$, in the 
usual sense (see for instance \cite[p. 182]{DiP}); accordingly it is uniquely determined by $f$, so that 
its 1-sum $f_{\rho_0}(x)$ is also uniquely determined by $f$. Applying the induction procedure to 
$f(x) - e^{ \rho_0 x} f_{\rho_0}(x)$ with $\Sigma\setminus\{\rho_0\}$ concludes the proof of Lemma \ref{lemasyunique}.
\end{proof}

\subsection{Notation and statement of Theorem \ref{theoprecis}} \label{subsecnotationsinf}

We consider a non-polynomial $E$-function $E(x)$ such that $E(0)= 0$, and write
$$
E(x)=\sum_{n=1}^\infty \frac{a_n}{n!}x^n.
$$
Its associated $G$-function is
$$
G(z)=\sum_{n=1}^\infty a_n z^n.
$$
We denote by $\calD$ a $G$-operator such that $\calFbar \calD E = 0$, where 
$\calFbar : \C[z,\frac{\dd}{\dd z}] \to \C[x,\frac{\dd}{\dd x}]$ is the Fourier transform 
of differential operators, i.e. the morphism of $\C$-algebras defined by 
$\calFbar(z) = \frac{\dd}{\dd x}$ and $\calFbar(\frac{\dd}{\dd z}) = -x$. Recall that such 
a $\calD $ exists because $E$ is annihilated by an $E$-operator, and any $E$-operator can be 
written as $\calFbar \calD$ for some $G$-operator $\calD$.

We let $g(z)=\frac{1}{z}G(\frac1z)$, so that $(\frac{\dd}{\dd z })^\delta \calD g = 0$ where 
$\delta$ is the degree of $\calD$ (i.e. the order of $\calFbar \calD$; see \cite{YA1}, p. 716). 
This function is the Laplace transform of $E(x)$: for $\Re(z)>C$, where $C>0$ is such that $\vert a_n\vert \ll C^{n}$, we have
$$
g(z)= \int_0^\infty E(x) e^{-xz} \dd x.
$$
From the definition of $g(z)$ and the assumption $E(0)=0$ we deduce that 
 $g(z)=\mathcal{O}(1/\vert z\vert ^2)$ as $z\to \infty$.

\medskip

We denote by 
$\nvsing$ the set of all finite singularities $\rho $ of $ \calD$; observe that 
$(\frac{\dd}{\dd z })^\delta \calD$ has the same singularities as $ \calD$. 
We also let 
$$\horsspec = \R\setminus\{\arg(\rho-\rho'), \rho,\rho'\in\nvsing, \rho\neq\rho'\}$$
where all the values modulo $2\pi$ of the argument of $\rho-\rho'$ are considered, so that $\horsspec+\pi = \horsspec$.

The directions $\theta\in\R\setminus(-\horsspec)$ (i.e., such that $(\rho-\rho')e^{i\theta}$ 
is real for some $\rho\neq\rho'$ in $\nvsing$) may be {\em anti-Stokes} (or {\em singular}, 
see for instance \cite[p. 79]{Loday}): when crossing such a direction, the renormalized sum of 
a formal solution at infinity of $\calD$ may change. In this paper we restrict to directions $\theta\in-\horsspec$.

For any $\rho\in\nvsing$ we denote by $\deltarho = \rho - e^{-i\theta}\R_+$ the half-line of 
angle $-\theta+\pi \bmod 2\pi$ starting at $\rho$. Since
$-\theta\in\horsspec$, no singularity $\rho'\neq\rho$ of $ \calD$ lies on $\deltarho $: these 
half-lines are pairwise disjoint. We shall work in the simply connected cut plane obtained from $\C$ 
by removing the union of these closed half-lines. We agree that for $\rho \in \nvsing$ and $z$ in the 
cut plane, $\arg(z-\rho)$ will be chosen in the open interval $(-\theta-\pi,-\theta+\pi)$. This enables 
one to define $\log(z-\rho)$ and $(z-\rho)^\al$ for any $\al \in \Q$.

Now let us fix $\rho\in\nvsing$. Combining theorems of Andr\'e, Chudnovski and Katz (see \cite[p. 719]{YA1}), 
there exist (non necessarily distinct) rational numbers $t_1^\rho, \ldots, t_{J(\rho)}^\rho$, with $J(\rho)\geq 1$, 
and $G$-functions $g_{j,k}^\rho$, for $1\leq j \leq J(\rho) $ and $0\leq k \leq K(\rho,j)$, such that a basis of local solutions of $(\frac{\dd}{\dd z })^\delta \calD$ around $\rho$ (in the above-mentioned cut plane) 
is given by the functions
 \begin{equation}\label{eqdeffjk}
f_{j,k}^\rho(z-\rho) = (z-\rho)^{t_j^\rho} \sum_{k'=0}^k g_{j,k-k'}^\rho(z-\rho) \frac{(\log(z-\rho))^{k'}}{k'!}
\end{equation}
for $1\leq j \leq J(\rho) $ and $0\leq k \leq K(\rho,j)$. Since $(\frac{\dd}{\dd z })^\delta \calD g= 0$ 
we can expand $g$ in this basis:
\begin{equation} \label{eqdefccg}
g(z) = \sum_{j=1}^{J(\rho)}\sum_{k=0}^{K(\rho,j)}\varpi_{j,k}^\rho f_{j,k}^\rho(z-\rho)
\end{equation}
with connexion constants $\varpi_{j,k}^\rho$; Theorem 2 of \cite{gvalues} yields $\varpi_{j,k}^\rho\in\GG$.

We denote by $\{u\} \in[0,1)$ the fractional part of a real number $u$, and agree that all 
derivatives of this or related functions taken at integers will be right-derivatives. We also denote
 by $\hada$ the Hadamard (coefficientwise) product of formal power series in $z$, and we let
$$y_{\al,i}(z) = \sum_{n=0}^\infty \frac{1}{i!} \frac{\dd^{i}}{\dd y^{i}}\Big(\frac{\Gamma(1-\{y\})}
{\Gamma(-y-n)}\Big)_{| y=\al } z^n \in\Q[[z]]$$
for $\al\in\Q$ and $i\in\N$. To compute the coefficients of $y_{\al,i}(z) $, we may restrict to 
values of $y$ with the same integer part as $\al$, denoted by $\peal$. Then 
\begin{equation} \label{eqconcretun}
\frac{\Gamma(1-\{y\})}{\Gamma(-y-n)} = \frac{\Gamma( - y +\peal +1)}{\Gamma(-y-n)} = 
\left\{ \begin{array}{l}
(-y-n)_{n+\peal + 1} \mbox{ if } n\geq - \peal \\ 
\\
\frac1{(-y+\peal+1)_{-n-\peal-1}}\mbox{ if } n \leq -1-\peal
\end{array}\right.
\end{equation}
is a rational function of $y$ with rational coefficients, so that $y_{\al,i}(z) \in \Q[[z]]$. 
Even though this won't be used in the present paper, we mention that $y_{\al,i}(z)$ is an arithmetic Gevrey series of order~$1$ (see \cite{YA1}); in particular it is divergent for any $z\neq 0$ (unless it is a polynomial, namely if $i=0$ and $\alpha\in\Z$).

 Finally, we define 
$$\eta_{j,k}^\rho(1/x) = \sum_{m=0}^k (y_{t_j^\rho,m}\hada g_{j,k-m}^\rho)(1/x) \in \Qbar[[1/x]]$$
for any $1\leq j \leq J(\rho)$ and $0\leq k \leq K(j,\rho)$; this is also an arithmetic Gevrey series of order~$1$. 
It is not difficult to see that 
$\eta_{j,k}^\rho(1/x) = 0$ if $f_{j,k}^\rho(z-\rho)$ is holomorphic at $\rho$. Indeed in this 
case $k=0$ and $t_j^\rho\in\Z$; if $t_j^\rho\geq 0$ then $y_{t_j^\rho,0}$ is identically zero, 
and if $t_j^\rho\leq -1$ then $y_{t_j^\rho,0}$ is a polynomial in $z$ of degree $-1-t_j^\rho$ 
whereas $g_{j,0}^\rho$ has valuation at least $-t_j^\rho$.

The main result of this section is the following asymptotic expansion, valid in the setting of 
Definition~\ref{defiasy} for $\theta\in-\horsspec$. It is at the heart of Theorem~\ref{thintrocce}; recall that we assume here $E(0)=0$, and that we let $\Gh=1/\Gamma$.

\begin{theo}\label{theoprecis} We have
$$E(x) \devasy \sum_{\rho\in\Sigma} e^{\rho x} \sum_{j=1}^{J(\rho)} \sum_{k=0}^{K(j,\rho)} \varpi_{j,k}^\rho x^{-t_j^\rho -1}
 \sum_{i = 0}^{ k} \Big( \sum_{\ell = 0} ^{k -i} \frac{(-1)^{\ell}}{\ell!} 
\Gh^{(\ell)}(1-\{t_j^\rho \}) \eta_{j, k-\ell-i}^\rho (1/x) \Big) \frac{ (\log(1/x))^{i} }{i!}.$$
\end{theo}
We observe that the coefficients are naturally expressed in terms of $\Gh^{(\ell)}$. 
Let us write Theorem \ref{theoprecis} in a slightly different way.
For $t\in\Q $ and $s\in\N$, let 
$$\lambda_{t,s}(1/x) = \sum_{\nu=0}^s \frac{(-1)^{s-\nu}}{(s-\nu)!}
\Gh^{(s-\nu)}(1-\{t \}) \frac{ (\log(1/x))^\nu}{\nu!}.$$
 In particular, $\lambda_{t,0}(1/x) = \Gh(1-\{t\}) $ and $\lambda_{t,1}(1/x) 
= \Gh(1-\{t\}) \log(1/x) - \Gh^{(1)}(1-\{t\}) $;
 for $t\in\Z$ we have $\lambda_{t,1}(1/x) = \log(1/x) - \gamma$. 

Then Theorem \ref{theoprecis} reads (by letting $s = i+\ell$):
\begin{equation} \label{eqgammalog}
E(x) \devasy \sum_{\rho\in\nvsing} e^{\rho x} \sum_{j=1}^{J(\rho)} \sum_{k=0}^{K(j,\rho)} 
\varpi_{j,k}^\rho x^{-t_j^\rho-1}\sum_{s=0}^k \lambda_{t_j^\rho,s}(1/x) \eta_{j, k-s}^\rho(1/x) .
\end{equation}

Here we see that the derivatives of $1/\Gamma$ do not appear in an arbitrary way, but 
always through these sums $\lambda_{t,s}(1/x) $. In particular $\gamma$ appears through $\lambda_{t,1}(1/x) = \log(1/x) - \gamma$, as mentioned in the introduction. 

\medskip

In the asymptotic expansion of Theorem \ref{theoprecis}, and in \eqref{eqgammalog}, the 
singularities $\rho\in\nvsing$ at which $g(z)$ is holomorphic have a zero contribution because 
for any $(j,k)$, either $\varpi_{j,k}^\rho=0$ or $f_{j,k}^\rho(z-\rho)$ is holomorphic at $\rho$ 
(and in the latter case, $k=0$ and $\eta_{j,0}^\rho(1/x) = 0$, as mentioned before the statement 
of Theorem \ref{theoprecis}). Moreover, as the proof shows (see \S \ref{subsecdemtheoprecis}), 
it is not really necessary to assume that the functions $f_{j,k}^\rho(z-\rho)$ form a basis of 
local solutions of $(\frac{\dd}{\dd z })^\delta \calD$ around $\rho$. Instead, it is enough to 
consider rational numbers $t_j^\rho$ and $G$-functions $g_{j,k}^\rho$ such that all singularities of 
 $g_{j,k}^\rho(z-\rho)$ belong to $\nvsing$ and, upon defining $f_{j,k}^\rho$ by Eq. \eqref{eqdeffjk}, 
Eq. \eqref{eqdefccg} holds with some complex numbers $\varpi_{j,k}^\rho$. In this way, to compute 
the asymptotic expansion of $E(x) $ it is not necessary to determine $\calD$ explicitly. 
The finite set $\nvsing$ is used simply to control 
the singularities of the functions which appear, and prevent $\theta$ from being a possibly singular direction. This remark makes it easier to apply Theorem \ref{theoprecis} to specific $E$-functions, for instance to obtain the expansions \eqref{eqex1} and \eqref{eqex2} used in \S \ref{sec:eapprox}. 

\subsection{Proof of Theorem \ref{theoprecis}} \label{subsecdemtheoprecis}

We fix an oriented line $\droite$ such that the angle between $\R_+$ and $\droite$ is equal 
to $-\theta+\frac{\pi}2 \bmod 2\pi$, and all singularities of $ \calD$ lie on the left of $\droite$. Let
 $R > 0$ be sufficiently large (in terms of $\droite$ and $\nvsing$). Then the circle $\cercle$ 
centered at 0 of radius $R$ intersects $\droite $ at two distinct points $a$ and $b$, with 
$\arg(b-a) = -\theta+\frac{\pi}2 \bmod 2\pi$, and 
\begin{equation}\label{eq:EG}
E(x)=\lim_{R\to\infty} \frac1{2i\pi}\int_a ^b g(z) e^{zx} \dd z 
\end{equation}
where the integral is taken along the line segment $ab$ contained in $\droite$.

\begin{figure} 
\centering
\includegraphics[width=0.5\textwidth ,trim = 0 320 0 80, clip=true]{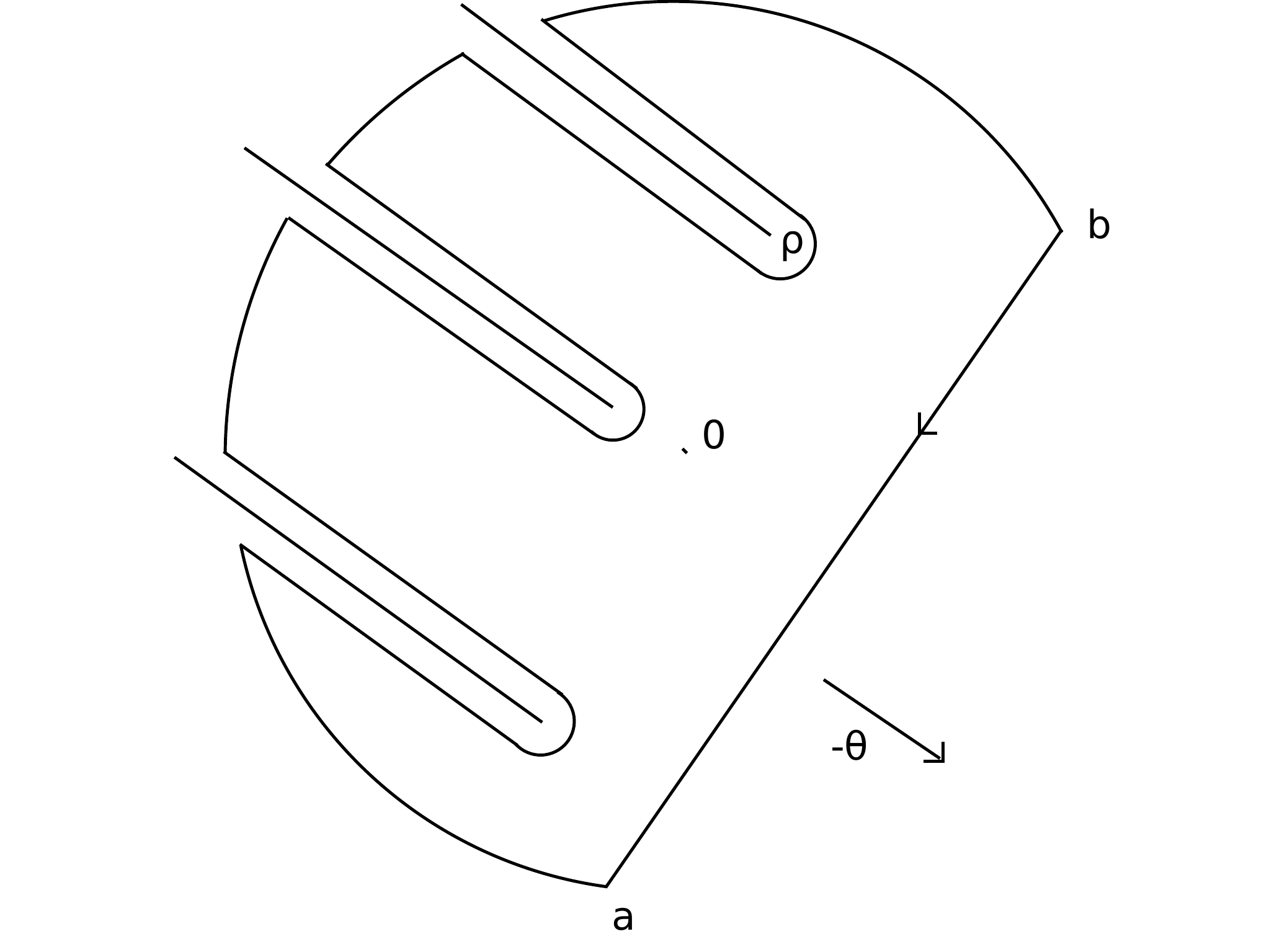}
\caption{The contour $\Gamma_R$} \label{fig1}
\end{figure}

 For any $\rho \in \nvsing$ the circle $\cercle $ intersects $\deltarho$ at one point 
$\pointrho = \rho - \Arho e^{ -i\theta}$, with $\Arho >0$, which corresponds to two points 
at the border of the cut plane, namely $\rho + \Arho e^{i(-\theta\pm\pi)}$ with values $-\theta\pm\pi$ 
of the argument. We consider the following path $\cheminrhoR$: a straight line from 
$\rho + \Arho e^{i(-\theta-\pi)}$ to $\rho$ (on one bank of the cut plane), then a circle around $\rho$ 
with essentially zero radius and $\arg(z-\rho)$ going up from $-\theta-\pi$ to $-\theta+\pi$, and 
finally a straight line from $\rho$ to $\rho + \Arho e^{i(-\theta+\pi)}$ on the other bank of the cut 
plane. We denote by $\groschemin$ the closed loop obtained by concatenation of the line segment $ba$, 
the arc $a\pointrhoun$ of the circle $\cercle$, the path $\cheminrhounR$, the arc $\pointrhoun\pointrhode$, 
the path $\cheminrhodeR$, \ldots, and the arc $\pointrhop b $ (where $\rho_1,\ldots,\rho_p $ are the 
distinct elements of $\nvsing$, ordered so that 
$\pointrhoun$, $\pointrhode$, \ldots, $\pointrhop$ are met successively
 when going along $\cercle$ from $a$ to $b$ in the negative direction); see Figure~\ref{fig1}. We refer to \cite[pp. 183--192]{DiP} for a similar computation.

We observe that
$$
\frac{1}{2i\pi} \int_{\groschemin} g(z) e^{zx} \dd z = 0
$$
for any $x\in\C$, because $\groschemin$ is a closed simple curve inside which the integrand has no singularity. 

From now on, we assume that $ \theta-\frac{\pi}2 < \arg(x) < \theta+\frac{\pi}2$.
As $R\to\infty$, the integral of $ g(z) e^{zx} $ over the line segment $ba$ tends to 
 $-E(x)$, using Eq. \eqref{eq:EG}. Moreover, as $z$ describes $\cheminrhoR$ (except maybe 
in a bounded neighborhood of $\rho$) we have $\Re(zx)<0$ and 
 $g(z)=\mathcal{O}(1/\vert z^2\vert)$, so that letting $R\to\infty$ one obtains (as in \cite{DiP})
 \begin{equation}\label{eq:Esectangul}
E(x)=\sum_{\rho\in\nvsing} \frac{1}{2i\pi}\int_{\cheminrhoinf} g(z) e^{zx} \dd z,
\end{equation}
where $ \cheminrhoinf$ is the 
 extension of $\cheminrhoR$ as $R\to \infty$.

Plugging Eq. \eqref{eqdefccg} into Eq. \eqref{eq:Esectangul} yields
 \begin{equation}\label{eq54bis}
 E(x) = \sum_{\rho\in\nvsing} \sum_{j=1}^{J(\rho)} \sum_{k=0}^{K(j,\rho)} \varpi_{j,k}^\rho \frac{1}{2i\pi}\int_{\cheminrhoinf} 
 f_{j,k}^\rho(z-\rho)e^{zx} \dd z .
\end{equation}
To study the integrals on the right hand side we shall prove the following general claim.
 {\em Let $\rho\in\nvsing$, and $\fct$ be a $G$-function such that $\fct(z-\rho)$ is 
 holomorphic on the cut plane. For any $\al\in\Q $ and any $k \in \N$, let 
$$\fctalk ( z-\rho) = \fct(z-\rho) (z-\rho)^{\al} \frac{ (\log(z-\rho))^{k }}{k!} .$$
Then
$$ \frac{1}{2i\pi}\int_{\cheminrhoinf} \fctalk ( z-\rho) e^{zx} \dd z $$
admits the following asymptotic expansion in a large sector bisected by $\theta$ (with $\Gh:=1/\Gamma$):
 $$e^{\rho x} x^{-\al -1} \sum_{\ell= 0} ^k \frac{(-1)^{\ell}}{\ell!} \Gh^{(\ell)}(1-\{\al \}) 
 \sum_{i = 0}^{k-\ell}\Big( y_{\al,k-\ell-i}\hada \fct \Big)(1/x)
 \frac{ (\log(1/x))^{i} }{i!}.$$}

To prove this claim, we first observe that
$$ \int_{\cheminrhoinf} \fctalk ( z-\rho) e^{zx} \dd z = \frac1{k!} \frac{\partial^{k}}{\partial \al^{k}} 
\Big[ \int_{\cheminrhoinf} \fctalzero ( z-\rho) e^{zx} \dd z\Big] 
$$
where the $k$-th derivative is taken at $\al$; this relation enables us to deduce the general case 
from the special case $k=0$ considered in \cite{DiP}. We write also
$$ \fct(z-\rho) = \sum_{n=0}^\infty c_{n} (z-\rho)^n. $$
Following \cite[pp. 185-191]{DiP}, given $\eps>0$ we obtain $R,C,\kappa>0$ such that, for any 
$n\geq 1$ and any $x$ with $|x| \geq R$ and $ \theta-\frac{\pi}2+\eps < \arg(x) < \theta+\frac{\pi}2-\eps$, we have
$$\Big| \frac{x^{-\al -n-1}}{\Gamma(-\al -n)} - \frac1{2i\pi}e^{-\rho x} 
\int_{\cheminrhoinf} (z-\rho )^{\al +n}e^{zx} \dd z \Big| \leq C^n n! |x|^{-\al-n-1}e^{-\kappa |x| \sin(\eps)}.$$
Then following the proof of \cite[pp. 191-192]{DiP} and using the fact that $\limsup |c_n|^{1/n}<\infty$ because $\fct$ is a $G$-function, for any
 $\eps>0$ we obtain $R,B, C >0$ such that, 
for any $N\geq 1$ and any $x$ with $|x| \geq R$ and $ \theta-\frac{\pi}2+\eps < \arg(x) < \theta+\frac{\pi}2-\eps$, we have
\begin{equation} \label{eqdevasyun}
\Big| e^{-\rho x} \frac{1}{2i\pi}\int_{\cheminrhoinf} \fctalk(z-\rho) e^{zx} \dd z - 
 \sum_{n=0}^{N-1} \frac{ c_n }{k!} \frac{\partial^{k }}{\partial \al^{k }} 
\Big[ \frac{x^{-\al-n-1}}{\Gamma(-\al-n)}\Big] \Big| \leq C^N N! |x|^{B-N}.
\end{equation}
Now observe that $\horsspec$ is a union of open intervals, so that $\theta$ can be made slightly 
larger or slightly smaller while remaining in the same open interval. In this process, the cut 
plane changes but the left handside of \eqref{eqdevasyun} remains the same (by the residue theorem, 
since $\fct(z-\rho)$ is holomorphic on the cut plane). The asymptotic expansion \eqref{eqdevasyun} 
remains valid as $|x| \to \infty$ in the new sector $ \theta-\frac{\pi}2+\eps < \arg(x) < \theta+\frac{\pi}2-\eps$, 
so that finally it is valid in a large sector 
$ \theta-\frac{\pi}2-\eps \leq \arg(x) \leq \theta+\frac{\pi}2+\eps$ for some $\eps>0$.

Now Leibniz' formula yields the following equality between functions of $\al$:
\begin{align*}
 \Big( \frac{x^{-\al-n-1}}{\Gamma(-\al-n)}\Big)^{(k )} &= \sum_{\ell = 0} ^k \sum_{i = 0}^{k-\ell} 
\frac{k !}{ \ell! i! (k-\ell-i)!} \big(\Gh(1-\{\al\})\big)^{(\ell)} 
 \Big( \frac{\Gamma(1-\{\al\}) }{\Gamma(-\al-n)}\Big)^{(k-\ell-i)}
\\
& \hspace{3cm} \times (\log(1/x))^{i} x^{-\al-n-1}
\\
&= \sum_{\ell = 0} ^k \frac{k!}{\ell!} \big(\Gh(1-\{\al\}\big)^{(\ell)} 
 \sum_{i = 0}^{k-\ell}\Big( y_{\al,k-\ell-i}\hada z^{n}\Big)(1/x)
 x^{-\al -1} \frac{ (\log(1/x))^{i} }{i!}
\end{align*}
 so that 
 $$ 
 \sum_{n=0}^{\infty} \frac{ c_n }{k!} \Big( \frac{x^{-\al-n-1}}{\Gamma(-\al-n)}\Big)^{(k)} = 
 \sum_{\ell = 0} ^k \frac1{\ell!} \big(\Gh(1-\{\al\})\big)^{(\ell)} 
 \sum_{i = 0}^{k-\ell}\Big( y_{\al,k-\ell-i}\hada \fct \Big)(1/x)
 x^{-\al -1} \frac{ (\log(1/x))^{i} }{i!}.
$$
 Using \eqref{eqdevasyun} this concludes the proof of the claim. 
 
 \medskip

Now we apply the claim to the $G$-functions $g_{j,k}^\rho$, since all singularities of
 $g_{j,k}^\rho(z-\rho)$ are singularities of $(\frac{\dd}{\dd z })^\delta \calD$ and
 therefore belong to $\nvsing$. Combining this result with Eqns. \eqref{eqdeffjk} and \eqref{eq54bis} yields:
\begin{align*}
E(x) &=
\sum_{\rho,j,k,k'}\varpi_{j,k}^\rho \frac{1}{2i\pi}\int_{\cheminrhoinf}
g_{j,k-k'}^\rho (z-\rho)(z-\rho)^{t_j^\rho} \frac{(\log(z-\rho))^{k'}}{k'!} e^{zx} \dd z \\
&\devasy\sum_{\rho,j,k,k'}\varpi_{j,k}^\rho e^{\rho x} x^{-t_j^\rho-1} \sum_{\ell = 0} ^{k'}
 \frac{(-1)^\ell}{\ell!} \Gh^{(\ell)}(1- \{t_j^\rho \})
 \sum_{i = 0}^{k'-\ell}\Big( y_{t_j^\rho, k'-\ell-i}\hada g_{j,k-k'}^\rho \Big)(1/x)
 \frac{ (\log(1/x))^{i} }{i!}
\\
&=\sum_{\rho,j,k}\varpi_{j,k}^\rho e^{\rho x} x^{-t_j^\rho-1} \sum_{\ell = 0} ^{k } 
\frac{(-1)^\ell}{\ell!} \Gh^{(\ell)}(1- \{t_j^\rho \})
 \sum_{i = 0}^{ k-\ell} \eta_{j, k-\ell-i}(1/x) \frac{ (\log(1/x))^{i} }{i!}.
\end{align*}

This concludes the proof of Theorem \ref{theoprecis}.

\subsection{Proof of Theorem \ref{thintrocce}} \label{subsec34}

To begin with, let us prove assertions $(ii)$ and $(iii)$. Adding the constant term $E(0) \in\Qbar\subset\GG$ to \eqref{eqdevintro} if necessary, we may assume that $E(0)=0$. Then Theorem \ref{theoprecis} applies; moreover, 
 in the setting of \S \ref{subsecnotationsinf} we may assume that the rational
 numbers $t_j^\rho$ have different integer parts as soon as they are distinct. Then 
letting $S$ denote the set of all $t_j^\rho+1$, for $\rho\in\nvsing$ and $1\leq j \leq J(\rho)$, 
and denoting by $T $ the set of non-negative integers less than or equal to $\max_{j,\rho} K(j,\rho)$, 
the asymptotic expansion of Theorem \ref{theoprecis} is exactly \eqref{eqdevintro} with coefficients
\begin{eqnarray*}
c_{\rho, \al,i,n} &=& \sum_{\stackrel{1\leq j \leq J(\rho)}{\mbox{{\tiny with }} \al = t_j^\rho+1}} \sum_{k=i}^{K(j,\rho)} \varpi_{j,k}^\rho 
 \sum_{\ell=0}^{k-i} \frac{(-1)^{\ell}}{\ell!} \Gh^{(\ell)}(1-\{\al \}) \\
&& \hspace{2cm} \sum_{m=0}^{k-\ell-i} \frac{1}{m!} \frac{\dd^{m}}{\dd y^{m}}
\Big(\frac{\Gamma(1-\{y\})}{\Gamma(-y-n)}\Big)_{| y=\al -1}
g_{j,k-\ell-i-m,n}^\rho 
\end{eqnarray*}
where $g_{j,k-\ell-i-m}^\rho(z-\rho) = \sum_{n=0}^\infty g_{j,k-\ell-i-m,n}^\rho (z-\rho)^n$. 
Now the coefficients $g_{j,k-\ell-i-m,n}^\rho$ are algebraic because $g_{j,k-\ell-i-m}^\rho$ is a 
$G$-function, and $\frac{\dd^{m}}{\dd y^{m}}\Big(\frac{\Gamma(1-\{y\})}{\Gamma(-y-n)}\Big)_{| y=\al -1}
$ is a rational number. Since $\varpi_{j,k}^\rho \in\GG$ and $\Qbar\subset\GG$, the coefficient $c_{\rho, \al,i,n} $ is a 
$\GG$-linear combination of derivatives of $\Gh = 1/\Gamma$ taken at the rational point $1-\{\al\}$. 
By the complements formula, $\Gh(z) = \frac{\sin(\pi z)}{\pi}\Gamma(1-z)$: applying Leibniz' formula we see 
that $\Gh^{(k)}(z)$ is a $\GG$-linear combination of derivatives of $\Gamma$ at $1-z$ up to order $k$, 
provided $z\in\Q\setminus\Z$ (using the fact \cite{gvalues} that $\GG$ contains $\pi$, $1/\pi$, 
and the algebraic numbers $\sin(\pi z)$ and $\cos(\pi z)$). When $z=1$, we use the identity (at $x=0$)
 $$\Gamma(x+1) = \exp\Big(-\gamma x + \sum_{k=2}^\infty \frac{(-1)^k \zeta(k)}{k}x^k\Big)$$
 (see \cite[p. 3, Theorem 1.1.2]{Andrews}) and the properties of Bell polynomials (see for 
instance \cite[Chap. III, \S 3]{Comtet1}). Since $\zeta(k) \in\GG$ for any $k\geq 2$ (because 
polylogarithms are $G$-functions), it follows that both $\Gamma^{(k)}(1)$ and $\Gh^{(k)}(1)$ 
are polynomials of degree $k$ in Euler's constant $\gamma$, with coefficients in $\GG$; moreover 
the leading coefficients of these polynomials are rational numbers. This implies that
 $\Gh^{(k)}(1)$ is a $\GG$-linear combination of derivatives of $\Gamma$ at 1 up to order $k$, and 
concludes the proof that all coefficients
 $c_{\rho, \al,i,n} $ in the expansion \eqref{eqdevintro} provided by Theorem~\ref{theoprecis} belong to ${\bf S}$.

To prove $(iii)$, we fix $\rho $ and $\alpha$ and denote by $K$ the maximal value of $K(j,\rho)$ among integers $j$ such that $\alpha = t_j^\rho+1$. Then 
$$c_{\rho, \al,i,n} = \sum_{\ell = 0}^{K-i} \frac{(-1)^\ell}{\ell!} \Gh^{(\ell)}(1-\{\alpha\}) g'_{\ell+i,n}$$
where
$$g'_{\lambda,n} = \sum_{j} \sum_{k=\lambda}^{K(j,\rho)} \varpi_{j,k}^\rho \sum_{m=0}^{k-\lambda} \frac1{m!} \frac{\dd^{m}}{\dd y^{m}}
\Big(\frac{\Gamma(1-\{y\})}{\Gamma(-y-n)}\Big)_{| y=\al -1} g_{j,k-\lambda-m,n}^\rho \in \GG;$$
here $0 \leq \lambda \leq K$ and the first sum is on $j \in \{1,\ldots,J(\rho)\}$ such that $\alpha = t_j^\rho+1$ and $K(j,\rho) \geq \lambda$. If $n$ is fixed and $g'_{\lambda,n}\neq 0$ for some $\lambda$, then denoting by $\lambda_0$ the largest such integer $\lambda$ we have $c_{\rho, \al,\lambda_0,n}\in \Gh(1-\{\alpha\})\cdot \GG\setminus\{0\} = \Gamma(\alpha)\cdot \GG\setminus\{0\}  $ and assertion $(iii)$ follows.

To prove $(i)$ and $(iv)$, we first observe that if $F(z)$ is given by  \eqref{eq:base0} with  algebraic  coefficients $\phiphi$, the asymptotic expansions of $F_j(z)$ we have just obtained can be multiplied by $\phiphi z^s \log(z)^k$ and summed up, thereby proving $(ii)$ for $F(z)$. To deduce $(i)$ from $(ii)$ for any solution $F(z)$ of an $E$-operator $L$, we recall that any   formal solution 
$f$ of $L$ at $\infty$   can be   written as \eqref{eqdevintro} with complex coefficients $c_{\rho, \al,i,n}(f) $, and  denote by $\listecoeff(f)$ the family of all these coefficients.
 The linear map $\listecoeff$ is injective, so that there exists a finite subset $\partiexi$ of the 
set of indices $(\rho, \al,i,n)$ such that $\petiteliste : f \mapsto ( c_{\rho, \al,i,n}(f))_{(\rho, \al,i,n) \in \partiexi}$
 is a bijective linear map. Denoting by $F_\theta$ the asymptotic expansion of $F(x)$ in a large sector
 bisected by $\theta$, we have
 $$\petiteliste(F_\theta) =\cce_{1,\theta} \petiteliste(H_1) + \ldots + \cce_{\mu,\theta} \petiteliste(H_\mu) $$
with the notation of \eqref{eqintro17}.
 Now $ \petiteliste(H_1) $, \ldots, $ \petiteliste(H_\mu)$ are linearly independent elements of
 $\Qbar^{\partiexi}$ and $\cce_{1,\theta} , \ldots, \cce_{\mu,\theta}$ can be obtained by Cramer's rule,
 so that they are linear combinations of the components of $\petiteliste(F_\theta)$ with coefficients 
in $\Qbar\subset\GG$: using $(ii)$ this concludes the proof of $(i)$.

\section{Asymptotics of the coefficients of $A(z)\cdot E\big(B(z)\big)$}\label{sec:asympPn}

In this section we deduce from Theorem \ref{thintrocce} the following result, of independent interest, which is the main step in the proof of Theorem \ref{theo:eapprox} (see \S \ref{preuvesubset2}).

\begin{theo} \label{theoasypn}
 Let $E(z)$ be an $E$-function, and $A(z), B(z) \in \Qbar[[z]]$ be algebraic functions; assume that $P(z) = A(z)\cdot E\big(B(z)\big) = \sum_{n=0}^{\infty} P_n z^n$ is not a polynomial. Then either
\begin{equation} \label{eqtheoaebun}
P_n = \frac{(2\pi)^{(1-d)/(2d)}}{n!^{1/d}}q^n n^{-u-1} (\log n)^v \Big( \sum_{\te} \Gamma(-u_\te) g_\te e^{in\te} + o(1)\Big)
\end{equation}
or
\begin{equation} \label{eqtheoaebde}
P_n = q^n e^{\sum_{\ell=1}^{d-1}\capa_\ell n^{\ell/d}} n^{-u-1} (\log n)^v \Big( \sum_{\te_1,\ldots,\te_d} \omega_{\te_1,\ldots,\te_d} e^{\sum_{\ell=1}^{d} i \te_\ell n^{\ell/d}} + o(1)\Big)
\end{equation}
where $q \in \Qbar$, $u \in \Q$, $u_\te\in \Q\setminus\N$, $d,v\in\N$, $d\geq 1$, $q > 0$, $g_\te\in \GG\setminus\{0\}$,
$\capa_1,\ldots,\capa_{d-1}\in\R$, $\te, \te_1,\ldots,\te_d\in [-\pi,\pi)$, the sums on $\te$ and $\te_1,\ldots,\te_d$ are finite and non-empty, and
\begin{equation} \label{eq43}
\left\{ \begin{array}{l}
 \omega_{\te_1,\ldots,\te_d} = \frac{\xi}{\Gamma(-u)} \mbox{ with } \xi\in ( \E \cup (\Gamma(\Q)\cdot \GG)) \setminus\{0\} \\
 \hspace{5cm}\mbox{ if } v=\capa_1= \ldots = \capa_{d-1}= \te_1=\ldots = \te_{d-1}=0,\\
 \omega_{\te_1,\ldots,\te_d} \in \Gamma(\Q) \cdot \exp( \Qbar)  \cdot \GG \setminus\{0\} \mbox{ otherwise.}
 \end{array}\right.
 \end{equation}
\end{theo}

As in the introduction, in \eqref{eq43} we let $\Gamma(-u) = 1$ if $u\in\N$.
In the special case where 
$$P(z) = (1-z)^{\alpha} \exp\Big(\sum_{i=1}^k \frac{b_i}{(1-z)^{\alpha_i}}\Big)$$
with $\alpha,\alpha_1,\ldots,\alpha_k\in\Q$, $b_1,\ldots,b_k\in\Qbar$, $\alpha_1 > 0$ and $b_1\neq 0$, 
Theorem~\ref{theoasypn} is consistent with Wright's asymptotic formulas \cite{wright2} for $P_n$.

\bigskip

We shall now prove Theorem \ref{theoasypn}; we distinguish between two cases (see \S \ref{subsecentiere} and \ref{subsecnonentiere}), which lead to Eqns. \eqref{eqtheoaebun} and \eqref{eqtheoaebde} respectively. This distinction, based on the growth of $P_n$, is different from the one mentioned in the introduction (namely whether $E(z)$ plays a role as $z\to z_0\in\C$ or as $z\to\infty$, providing elements of $\E$ or $\Gamma(\Q)\cdot\GG$ respectively). We start with the following consequence of Theorem \ref{thintrocce}, which is useful to study $E(z)$ as $z\to\infty$, in both \S \ref{subsecentiere} and \S \ref{subsecdiff}.

\begin{lem} \label{lemEsomme}
For any $E$-function $E(z)$ there exist $K \geq 1$, $u_1,\ldots,u_K\in\Q$, $v_1,\ldots,v_K\in\N$, and pairwise distinct $\alpha_1,\ldots,\alpha_K\in\Qbar$ such that 
\begin{equation} \label{eqdevEnv}
E(z) = \sum_{k=1}^K \omega_{k} e^{\alpha_k z} z^{u_k} \log(z)^{v_k} (1+o(1))
\end{equation}
as $|z|\to\infty$, uniformly with respect to $\arg(z)$, where $\omega_k \in \Gamma(-u_k) \cdot \GG\setminus\{0\} $ with  $\Gamma(-u_k)=1$ if $u_k\in\N$.
\end{lem}

In Eq. \eqref{eqdevEnv} we assume that a determination of $\log z$ is chosen in terms of $k$, with a cut in a direction where the term corresponding to $k$ is very small with respect to another one (except if $K=1$, but in this case the proof yields $v_1=0$ and $u_1\in\Z$). 

\begin{proof}
For any $\alpha\in\C$, let $I_\alpha$ denote the set of all directions $\theta\in\R/2\pi\Z$ such that $E(z)$ has an asymptotic expansion \eqref{eqdevintro} in a large sector bisected by $\theta$, with $\Sigma$ having the least possible cardinality, $\alpha\in\Sigma$, and $\Re(\alpha' e^{i\theta}) \leq \Re(\alpha e^{i\theta})$ for any $\alpha'\in\Sigma$. This implies that in the direction $\theta$, the growth of $E(z)$ is comparable to that of $e^{\alpha z}$. Then $I_\alpha$ is either empty or of the form   $[R_\alpha,S_\alpha] \bmod 2\pi$ with $R_\alpha\leq S_\alpha$. We denote by $\Sigma_0$ the set of all $\alpha\in\C$ such that $I_\alpha\neq\emptyset$; then $\Sigma_0$ is a subset of the finite set $\nvsing \subset\Qbar$ constructed in \S \ref{subsecnotationsinf}, so that $\Sigma_0$ is finite: we denote by $\alpha_1,\ldots,\alpha_K$ its elements, with $K\geq 1$.

If $K=1$ then $I_{\alpha_1} = \R/2\pi\Z$ and the asymptotic expansion \eqref{eqdevintro} is the same in any direction: $e^{-\alpha_1z} E(z)$ has (at most) a pole at $\infty$, and Lemma \ref{lemEsomme} holds with $u_1\in\Z$, $v_1=0$, and $\omega_1\in\GG$ (using Theorem \ref{thintrocce}).

Let us assume now that $K\geq 2$. Then $S_{\alpha_k} - R_{\alpha_k} \leq \pi$ for any $k$, so that $E(z)$ admits an asymptotic expansion \eqref{eqdevintro} in a large sector that contains all directions $\theta\in I_{\alpha_k}$. Among all terms corresponding to $e^{\alpha_k z}$ in this expansion, we denote the leading one by
\begin{equation}
\label{eqdefalk}
\omega_k e^{\alpha_k z} z^{u_k} (\log z)^{v_k}
\end{equation}
with $u_k\in\Q$, $v_k\in\N$, and $\omega_k \in  \Gamma(-u_k) \cdot \GG\setminus\{0\} $ (using assertion $(iii)$ of Theorem \ref{thintrocce}), where $\Gamma(-u_k)$ is understood as 1 if $u_k$ is a non-negative integer. These parameters are the one in \eqref{eqdevEnv}. To conclude the proof of Lemma \ref{lemEsomme}, 
we may assume that $\arg(z)$ remains in a small segment $I$, and consider the asymptotic expansion \eqref{eqdevintro} in a large sector containing $I$. Keeping only the dominant term corresponding to each $\alpha\in\Sigma$ in this expansion, we obtain
\begin{equation}
\label{eqdevEde}
E(z) = \sum_{\alpha\in\Sigma} \omega'_{\alpha} e^{\alpha z} z^{u'_\alpha} (\log z)^{v'_\alpha} (1+o(1)).
\end{equation}
To prove that \eqref{eqdevEde} is equivalent to \eqref{eqdevEnv} as $|z|\to\infty$ with $\arg(z)\in I$, we may remove from both equations all terms corresponding to values $\alpha_k$ (resp. $\alpha\in\Sigma$) such that $I_{\alpha_k} \cap I = \emptyset$ (resp. $I_{\alpha} \cap I = \emptyset$), since they fall into error terms. Now for any $\alpha = \alpha_k$ such that $I_{\alpha} \cap I \neq \emptyset$, $E(z)$ admits an asymptotic expansion in a large sector containing $I_{\alpha} \cup I$ (since $I_\alpha$ has length at most $\pi$, and the length of $I$ can be assumed to be sufficiently small in terms of $E$). Comparing the dominating exponential term of this expansion in a direction $\theta \in I_{\alpha} \cap I$ with the ones of \eqref{eqdefalk} and \eqref{eqdevEde}, we obtain $ \omega'_{\alpha} = \omega_k$, $u'_\alpha = u_k$, and $v'_\alpha = v_k$. This concludes the proof of Lemma \ref{lemEsomme}.
\end{proof}

\subsection{$P(z)$ is an entire function}\label{subsecentiere}

If $P(z)$ is an entire function then $A(z)$ and $B(z)$ are polynomials; we denote by $\delta \geq 0$ and $d \geq 1$ their degrees, and by $A_\delta$ and $B_d$ their leading coefficients. 
We shall estimate the growth of the Taylor coefficients of 
$P(z)$ by the saddle point method. For any circle $C_R$ of center $0$ and 
radius $R$, Lemma \ref{lemEsomme} yields
\begin{align}
P_n&=\frac{1}{2 i \pi} \int_{C_R} \frac{A(z) \cdot E(B(z))}{z^{n+1}} \dd z\nonumber \\
&=\frac{1}{2 i \pi} 
\sum_{k=1}^K \omega_k A_\delta B_d^{u_k} d^{v_k} \int_{C_R } e^{\alpha_k B(z)} \cdot z^{\delta+d u_k -n-1} (\log z)^{v_k} \cdot (1+o(1))\dd z \nonumber 
\end{align}
where the $o(1)$ is with respect to $R \to +\infty$ and is uniform in $n$; here $\log(z)$ is a fixed determination which depends on $k$ (see the remark after Lemma \ref{lemEsomme}). We have to distinguish between the cases $\alpha_k=0$ and $\alpha_k\neq 0$. In 
the former case, the integral 
$$
\frac{\omega_k}{2 i \pi}
 \int_{C_R } z^{\delta+d u_k -n-1} (\log z)^{v_k} \cdot (1+o(1))\dd z
$$
 tends to $0$ as $R\to +\infty$ (provided $n$ is sufficiently large) and there is no contribution coming from this case.
 
Now $E(z)$ is not a polynomial (otherwise $P(z)$ would be a polynomial too), so that if $\alpha_k = 0$ for some $k$ then $K \geq 2$: there is always at least one integer $k$ such that $\alpha_k\neq 0$. For any such $k$, 
the function 
$$
 e^{\alpha_k B(z)} z^{\delta+d u_k -n-1} (\log z)^{v_k}
$$
is smooth on $C_R $ (except on the cut of $\log z$) and the integral can be estimated
as $n\to \infty$ by finding the critical points of $\alpha_k B(z)-n \log(z)$, i.e. the solutions 
$z_{1,k}(n), \ldots, z_{d,k}(n)$ of $zB'(z)=n/\alpha_k$. As $n\to\infty$, we have 
$z_{j,k}(n)\equi (d B_{d} \alpha_k)^{-1/d} e^{2i \pi j/d} n^{1/d} \to \infty$, so that $\alpha_k B(z_{j,k}(n))\equi n/d$. 

Moreover, denoting by $\Delta_{j,k}(n)$ the second derivative of $\alpha_k B(z)-n \log(z)$ at $z=z_{j,k}(n)$, we see that asymptotically
$$
\Delta_{j,k}(n)=\alpha_k B''(z_{j,k}(n))+\frac{n}{z_{j,k}(n)^2} 
\equi d (d B_d \alpha_k)^{2/d}e^{-4i \pi j/d} n^{1-2/d}.
$$
Then the saddle point method yields:
$$P_n = \sum_{\alpha_k\neq 0} \omega'_k \sum_{j=0}^{d-1} \frac{1}{\sqrt {2\pi \Delta_{j,k}(n) }}
 e^{\alpha_k B(z_{j,k}(n))} z_{j,k}(n)^{\delta+d u_k-n-1} (\log z_{j,k}(n))^{v_k} (1+o(1)) $$
with $\omega'_k = \omega_k A_\delta B_d^{u_k} d^{v_k} \in \Qbaretoile \omega_k $. This relation yields
$$P_n = \sum_{\alpha_k\neq 0} \frac{\omega''_k}{\sqrt{2\pi}} n^{-n/d} (ed B_d\alpha_k)^{n/d} n^{\frac{\delta}{d}+ u_k -\frac12}(\log n)^{v_k} \Big( \sum_{j=0}^{d-1} e^{2i\pi jn/d} + o(1)\Big)$$
with $\omega''_k \in \Qbaretoile \omega_k $. Now let $\alti = \max(|\alpha_1|,\ldots,|\alpha_K|)$ and consider the set $\calK$ of all $k$ such that $|\alpha_k| = \alti$.
For each $k\in\calK$ we write $\alpha_k^{1/d} = \alti^{1/d} e^{i\te_k}$; then Stirling's formula yields
$$P_n = (2\pi)^{(1-d)/(2d)} n!^{-1/d} (dB_d\alti)^{n/d} \sum_{k\in\calK} \omega''_k n^{\frac{\delta}{d} +u_k -\frac12 + \frac1{2d}}(\log n)^{v_k} \sum_{j=0}^{d-1} e^{i (\te_k + \frac{2 \pi j}{d})n} (1 + o(1)).
$$
Keeping only the dominant terms provides Eq. \eqref{eqtheoaebun}.

\subsection{$P(z)$ is not an entire function} \label{subsecnonentiere}

Let us move now to the case where $P(z)$ is not entire. It has only a finite number of singularities of minimal modulus (equal to $q^{-1}$, say), and as usual the contributions of these singularities add up to determine the asymptotic behavior of $P_n$. Therefore, for simplicity we shall restrict in the proof to the case of a unique singularity $\rho$ of minimal modulus $q^{-1}$. We consider first two special cases, and then the most difficult one.

\subsubsection{$B(z)$ has a finite limit at $\rho$} \label{subsec421}

Let us assume that $B(z)$ admits a finite limit as $z\to \rho$, denoted by $B(\rho)$; $\rho$ can be a singularity of $B$ or not. In both cases, as $z \to \rho $ we have
$$B(z) = B(\rho) + \mathfrak{B} (z-\rho)^t (1+o(1))$$
with $t\in\Q$, $t\geq 0$, and $ \mathfrak{B}\in\Qbaretoile$ (unless $B$ is a constant; in this case the proof is even easier). Now all Taylor coefficients of $E(z)$ at $B(\rho)$ belong to $\E$, so that 
$$E(B(z)) \equi \eta (z-\rho)^{t'}$$
as $z\to\rho$, with $t'\in\Q$, $t'\geq 0$, and $\eta\in\E\setminus\{0\}$. On the other hand, if $\rho$ is a singularity of the algebraic function $A(z)$ then its
 Puiseux expansion yields $s\in\Q\setminus\N$, $\mathfrak{A} \in\Qbaretoile$ and a polynomial $\widetilde A$ such that 
$$A(z) = \widetilde A(z-\rho) + \mathfrak{A} (z-\rho)^s (1+o(1))$$
as $z\to\rho$; if $\rho$ is not a singularity of $A$ we have the same expression with $s\in\N$ and $\widetilde A = 0$. In both cases we obtain finally $p \in \Q\setminus\N$, $ \mathfrak{P} \in\E \setminus\{0\}$ and a polynomial $\widetilde P$ such that 
$$P(z) = \widetilde P(z-\rho) + \mathfrak{P} (z-\rho)^p (1+o(1)).$$
Using standard transfer results (see \cite{Bible}, p. 393) this implies
$$P_n \equi \frac{(-\rho)^{-p} \mathfrak{P} }{\Gamma(-p)} \rho^{-n}n^{-p-1}.$$
Therefore the singularity contributes to \eqref{eqtheoaebde} through a term in which $v=\capa_1= \ldots = \capa_{d-1}= \te_1=\ldots = \te_{d-1}=0$ and $\rho^{-1} = q e^{i\te_d}$. 

\subsubsection{$E$ is a polynomial} \label{subsec422}

In this case, $P(z)$ is an algebraic function (and not a polynomial) so that
$$
P_n\equi \frac{\omega}{\Gamma(-s)}\cdot n^{-s-1} \rho^{-n}
$$
with $\omega \in \Qbaretoile \subset \E\setminus\{0\}$ and $s\in \Q\setminus\N$ determined by the Puiseux expansion of $P(z)$ around $\rho$ (using the same transfer result as above). Therefore each singularity $\rho = q^{-1} e^{-i\te_d}$ contributes to a term in \eqref{eqtheoaebde} with $v=\capa_1= \ldots = \capa_{d-1}= \te_1=\ldots = \te_{d-1}=0$.

\subsubsection{The main part of the proof} \label{subsecdiff}

Let us come now to the most difficult part of the proof, namely the contribution of a singularity $\rho$ at which $B(z)$ does not have a finite limit (in the case where $E(z)$ is not a polynomial). As above we assume (for simplicity) that $\rho$ is the unique singularity of $P(z)$ of minimal modulus $q^{-1}$. As $z\to\rho$, we have
\begin{equation} \label{eqasyab}
A(z) \equi \mathfrak{A} (z-\rho)^{t/s} \mbox{ and }
 B(z) \equi \mathfrak{B} (z-\rho)^{-\tau/\sigma}
\end{equation}
with $\mathfrak{A} , \mathfrak{B} \in \Qbaretoile$, $s,t,\sigma,\tau\in \mathbb Z$, $s,\sigma,\tau>0$, and $\gcd(s,t) = \gcd(\sigma,\tau) = 1$. For any circle $C_R$ of center $0$ and 
radius $R<\vert \rho\vert$, we have (using Lemma \ref{lemEsomme} as in \S \ref{subsecentiere})
\begin{equation} \label{eqpnde}
P_n =\frac{1}{2 i \pi}
\sum_{k=1}^K \omega_k \int_{C_R } \frac{e^{\alpha_k B(z)}}{z^{n+1}}
\cdot A(z)B(z)^{u_k}\log(B(z))^{v_k} \cdot (1+o(1))\dd z
\end{equation} 
where $o(1)$ is with respect to $R \to \vert \rho\vert$ and is uniform in $n$. 

If $\alpha_k=0$ for some $k$, then the corresponding term in \eqref{eqpnde} has to be treated in a specific way, since the main contribution may come from the error term $o(1)$. For this reason we observe that in Lemma \ref{lemEsomme}, the term corresponding to $\alpha_k=0$ can be replaced with any truncation of the asymptotic expansion of $E(z)$, namely with
$$\sum_{u = -U_0}^{U_1} \sum_{v=0}^{V} \omega_{u,v} z^{u/d} (\log z)^v + o(z^{-U_0/d})$$
where $d\geq 1$ and $U_0$ can be chosen arbitrarily large. Now the corresponding term in \eqref{eqpnde} becomes
\begin{equation} \label{contribexcep}
\frac1{2i\pi} \int_{C_R } \frac1{z^{n+1}} \Big( \sum_{u = -U_0}^{U_1} \sum_{v=0}^{V} \omega_{u,v}A(z) B(z)^{u/d} (\log B(z))^v + o(A(z) B(z)^{-U_0/d})\Big)dz.
\end{equation}
The point is that the function $ \omega_{u,v}A(z) B(z)^{u/d} (\log B(z))^v $ may be holomorphic at $z=\rho$, because $ \omega_{u,v}=0$ or because the singularities at $\rho$ of $A(z)$ and $B(z)^{u/d} (\log B(z))^v $ cancel out; in this case the corresponding integral over $C_R$ is $o(q'^n)$ for some $q'<q=|\rho|^{-1}$ so that it falls into error terms. If this happens for any $U_0$, for any $u$ and any $v$, then the term corresponding to $\alpha_k=0$ in \eqref{eqpnde} is $o(q^nn^{-U})$ for any $U>0$, so that it falls into the error term of the expression \eqref{eqtheoaebde} we are going to obtain for $P_n$. Otherwise we may consider the maximal pair $(u,v)$ (with respect to lexicographic order) for which this function is not holomorphic; then \eqref{contribexcep} is equal to 
$$
 \frac{ \omega'_{u,v}}{2 i \pi} \int_{C_R} \frac{(\rho-z)^{T}\log(\rho-z)^{v }}{z^{n+1}}
\cdot (1+o(1))\dd z 
$$
for some $T\in\Q$ and $ \omega'_{u,v} \in \Qbaretoile \omega_{u,v}\subset \Gamma(\Q) \cdot \GG$ (using assertion $(iii)$ of Theorem \ref{thintrocce}). We obtain finally (see \cite{Bible}, p. 387): 
$$
\begin{cases}
\displaystyle \frac{ \omega'_{u,v}}{\Gamma(-T)}\rho^{T-n} n^{-T-1}\log(n)^{v} (1+o(1)) \quad \textup{if}\; T\not\in\N,
\\
\\
\displaystyle\omega'_{u,v} \rho^{T-n} n^{-T-1} \log(n)^{v -1} (1+o(1)) \quad \textup{if}\; T\in\N \mbox{ (so that } v\geq 1\mbox{).}
\end{cases}
$$
This contribution can either fall into the error term of \eqref{eqtheoaebde}, or give a term with 
$ \capa_1= \ldots = \capa_{d-1}= \te_1=\ldots = \te_{d-1}=0$.

\bigskip

Let us now study the terms in \eqref{eqpnde} for which $ \alpha_k\neq 0$; since $E(z)$ is not a polynomial there is at least one such term. The function 
$$
\frac{e^{\alpha_k B(z)}}{z^{n+1}}\cdot A(z)B(z)^{u_k}\log(B(z))^{v_k}
$$
is smooth on $C_R $ (except on the cuts of $\log(B(z))$) and the integral can be estimated
as $n\to \infty$ by finding the critical points of $\alpha_k B(z)-n \log(z)$, i.e. the solutions 
of $zB'(z)=n/\alpha_k$. For large $n$, any critical point $z$ must be close to $\rho$ 
(since $zB'(z)$ is bounded away from $\rho$ for $|z| \leq |\rho|$). Now in a neighborhood of $z=\rho$ we have
$$
zB'(z) \equi -\frac{\rho \tau \mathfrak{B}}{\sigma} \cdot \frac{1}{(z-\rho)^{1+\tau/\sigma}}
$$
so that we have $\tau+\sigma$ critical points $z_{j,k}(n)$, for $j=0, \ldots, \sigma+\tau-1$, with 
$$
z_{j,k}(n) - \rho \equi e^{2 i \pi j\sigma/(\sigma+\tau)}\cdot \left(-\frac{\sigma n}{\rho \mathfrak{B} \tau \alpha_k}
\right)^{-\sigma/(\sigma+\tau)}.
$$ 
Using \eqref{eqasyab} and letting $\kappa=t/s \in \mathbb Q$ we deduce that
$$
A(z_{j,k}(n)) \equi \mathfrak{A} e^{2 i \pi j\sigma\kappa/(\sigma+\tau)}\cdot 
\left(-\frac{\sigma n}{\rho \mathfrak{B} \tau \alpha_k}\right)^{-\sigma\kappa /(\sigma+\tau)} \neq 0.
$$
Moreover we have
$$\alpha_k B(z_{j,k}(n)) \equi \frac{-\sigma}{\tau} ( z_{j,k}(n) - \rho) \alpha_k B'(z_{j,k}(n)) \equi
\frac{-\sigma n}{\rho \tau} ( z_{j,k}(n) - \rho) \equi \ddjk n^{\tau/(\sigma+\tau)}$$
with 
\begin{equation} \label{eqdefdk}
\ddjk = \Big( \alpha_k \mathfrak{B} e^{2i\pi j}\Big)^{\sigma/(\sigma+\tau)} \Big(\frac{-\sigma}{\rho \tau}\Big)^{ \tau / (\sigma+\tau)}\neq 0.
\end{equation}
To apply the saddle point method, we need to estimate the second derivative 
$\Delta_{j,k}(n)$ of $\alpha_kB(z)-n\log(z)$ at $z=z_{j,k}(n)$. We obtain 
$$
\Delta_{j,k}(n)=\alpha_kB''(z_{j,k}(n))+\frac{n}{z_{j,k}(n)^2} \equi \frac{\tau(\sigma+\tau)}{\sigma^2} 
(\alpha_k \mathfrak{B})^{-\sigma/(\sigma+\tau)} e^{-2i\pi j\frac{2\sigma+\tau}{\sigma+\tau}}\left(-\frac{\sigma}
{\rho \tau }\right)^{\frac{ 2\sigma+\tau}{ \sigma+ \tau}}
n^{\frac{ 2\sigma+\tau}{ \sigma+ \tau}}.
$$
Finally, 
$$
B(z_{j,k}(n))^{u_k}\equi (\ddjk / \alpha_k)^{u_k} n^{\tau u_k / (\sigma+\tau)}.
$$
This enables us to apply the saddle point method. This yields a non-empty subset $J_k$ of $\{0,\ldots,\sigma+\tau-1\}$ such that the term corresponding to $\alpha_k$ in \eqref{eqpnde} is equal to 
$$
 \sum_{j\in J_k} \frac{ \omega_k }{\sqrt{2\pi \Delta_{j,k}(n)}}
\frac{e^{\alpha_k B(z_{j,k}(n))}}{z_{j,k}(n)^{n +1}} A(z_{j,k}(n))B(z_{j,k}(n))^{u_k} \log(B(z_{j,k}(n)))^{v_k} (1+o(1)) .$$
Now for any pair $(j,k)$, $\alpha_k B(z_{j,k}(n))$ is an algebraic function of $n$ so that it can be expanded as follows as $n\to\infty$:
\begin{equation} \label{eqasydk}
\alpha_k B(z_{j,k}(n)) = \sum_{\ell=0}^{d'} \capa_{j,k,\ell} n^{\ell / d} + o(1)
\end{equation}
with $\capa_{j,k,\ell} \in \Qbar$, $0<d'<d$ and $d'/d = \frac{\tau}{\sigma+\tau}$, $\capa_{j,k,d'} = \ddjk \neq 0$. Increasing $d$ and $d'$ if necessary, we may assume that they are independent from $(j,k)$. We denote by $(\capa_{d'},\ldots,\capa_{1})$ the family $(\Re \capa_{j,k,d'},\ldots, \Re\capa_{j,k,1})$ which is maximal with respect to lexicographic order (as $j$ and $k$ vary with $ \alpha_k\neq 0$ and $j\in J_k$), i.e. 
for which the real part of \eqref{eqasydk} has maximal growth as $n\to\infty$. Among the set of pairs $(j,k) $ for which $\Re \capa_{j,k,1} = \capa_1$, \ldots, $\Re \capa_{j,k,d'} = \capa_{d'}$, we define $\calK $ to be the subset of those for which $(u_k,v_k)$ is maximal (with respect to lexicographic order), and let $(u,v)$ denote this maximal value. 
Then the total contribution to \eqref{eqpnde} of all terms with $\alpha_k \neq 0$ is equal to 
$$ \frac{n^{- \frac{\tau+2(1+\kappa)\sigma}{2\tau+2\sigma}}}{\sqrt{2\pi}} \rho^{-n} n^{\tau u / (\sigma+\tau)} \log(n)^{v } 
e^{\sum_{\ell=1}^{d'}\capa_\ell n^{\ell/d}} \Big( \sum_{(j,k)\in \calK} \widehat{\omega}_{j,k} e^{\capa_{j,k,0}}
e^{\sum_{\ell=1}^{d'} i \Im \capa_{j,k,\ell} n^{\ell/d}} +o(1)\Big) 
$$
with $\widehat{\omega}_{j,k} \in \Qbaretoile \omega_k$. Since $\capa_{d'} + i \Im \capa_{j,k,d'} = \ddjk\neq 0$, this concludes the proof of Theorem \ref{theoasypn}.

\section{Application to $E$-approximations}\label{sec:eapprox}

In this section we prove the results on $E$-approximations stated in the introduction, and discuss in \S \ref{ssec:extended} the generalization involving \eqref{eq:eapproxgen}.

\subsection{Examples of $E$-approximations}\label{ssec:example}

We start with an emblematic example. The diagonal Pad\'e approximants to $\exp(z)$ are given by
$Q_n(z)e^z- P_n(z)= \mathcal{O}(z^{2n+1})$ with 
$$
Q_n(z)=\sum_{k=0}^n (-1)^{n-k}\binom{2n-k}{n}\frac{z^k}{k!} \quad \textup{and} \quad P_n(z)=-Q_n(-z).
$$
It is easy to prove that, for any $z\in \mathbb C$ and any $x$ such that $\vert x\vert<1/4$, 
$$
\sum_{n=0}^\infty Q_n(z) x^k = \frac{e^{- z/2}}{\sqrt{1+4x}}e^{\frac z2\sqrt{1+4x}}. 
$$
This generating function can be written as 
$
\frac{e^{-z/2}}{\sqrt{1+4x}}f(z,x) + g(z,x),
$
where $f(z,x)$ and $g(z,x)$ are entire functions of $x$, and $f(z,-\frac 14)=-f(-z,-\frac14)\neq 0$. Hence, the asymptotic behavior of $Q_n(z)$ and $P_n(z)$ are given by  
$$
Q_n(z)\sim e^{-z/2}f\Big(z,-\frac14\Big) 4^n\binom{-1/2}{n}
\quad \textup{and} \quad
P_n(z)\equi e^{z/2}f\Big(z,-\frac14\Big) 4^n\binom{-1/2}{n}.
$$
It follows in particular that 
$$
\lim_{n\to +\infty} \frac{P_n(z)}{Q_n(z)} =e^z.
$$
This proves that for any $z\in \Qbar$, $e^z$ has $E$-approximations. Moreover, it is well-known that $n!P_n(1)$ and $n!Q_n(1)$ are respectively the numerator and denominator  of the $n$-th convergent of the continued fraction of the number $e$. In other words, the convergents of $e$ are $E$-approximations of $e$.

\bigskip

As mentioned in the introduction, any element of $\Frac\GG$ has $E$-approximations. To complete the proof of \eqref{eq:subset1}, let us prove this for any element of 
$\frac{ \E \cup \Gamma(\Q)}{ \E \cup \Gamma(\Q)} $ by constructing for any $\xi \in \E \cup \Gamma(\Q)$ a sequence $(P_n)$ as in Definition \ref{def:Eapprox} with  $\lim_{n\to\infty} P_n = \xi$. 

If $\xi = F(\alpha)$ where $\alpha\in\Qbar$ and 
$F(z)=\sum_{n\ge 0} \frac{a_n}{n!}z^n$ is an $E$-function, we define $P_n\in \Qbar$ by 
$$
\sum_{n=0}^{\infty} P_n z^n = \frac{1}{1-z} F(\alpha z).
$$
Then, trivially, 
$$
P_n=\sum_{k=0}^n \frac{a_k}{k!}\alpha^k \longrightarrow F(\al) = \xi.
$$

If $\xi = \Gamma(\alpha)$ with 
$\alpha \in \mathbb Q\setminus \mathbb Z_{\le 0}$, we consider the $E$-function 
$$
E_\alpha(z)=\sum_{n=0}^{\infty} \frac{z^n}{n!(n+\alpha)} 
$$
and define 
$P_n(\alpha)$ as announced in the introduction, by the series expansion (for $\vert z\vert <1$)
$$
\frac1{(1-z)^{\alpha+1}} E_\alpha\left(-\frac{z}{1-z}\right) = \sum_{n\ge 0} P_n(\alpha) z^n \in \mathbb Q  [[z]].
$$
Then 
$$
P_n(\alpha)=\sum_{k=0}^n \binom{n+\al}{k+\al}\frac{(-1)^k}{k!(k+\al)}
$$
(by direct manipulations) 
and, provided that $\alpha<1$,
$$
\lim_{n\to +\infty} P_n(\alpha) = \Gamma(\alpha) = \xi .
$$
To see this, we start from the asymptotic expansion
\begin{equation} \label{eqex1}
E_\alpha(-z)\devasy \frac{\Gamma(\alpha)}{z^{\alpha}} - e^{-z}\sum_{n=0}^{\infty} (-1)^n 
\frac{(1-\alpha)_n}{z^{n+1}}
\end{equation}
in a large sector bisected by any $\theta\in(-\frac\pi2, \frac\pi2)$, which is a special case of Theorem \ref{theoprecis} (proved directly in \cite{Michigan}, Proposition 1). Since 
$\exp\big(-\frac z{1-z}\big)=\mathcal{O}(1)$, as $z\to 1$, $\vert z\vert <1$, 
it follows that 
\begin{equation*}
\frac1{(1-z)^{\alpha+1}} E_\alpha\left(-\frac{z}{1-z}\right) =
\frac{\Gamma(\alpha)}{1-z} + \mathcal{O}\left(\frac1{\vert 1-z\vert^{\alpha}}\right)
\end{equation*}
for $z\to 1$, $\vert z\vert<1$. The result follows by 
standard transfer theorems since  $\alpha<1$; this example is of the type covered by \S \ref{subsecdiff} with $\alpha_1=0$.

From the differential equation $zy''(z)+(\al+1-z)y'(z)-\al y(z)=0$ satisfied by $E_{\alpha}(z)$, we easily get the differential equation satisfied by $\frac1{(1-z)^{\alpha+1}} E_\alpha\left(-\frac{z}{1-z}\right)$: 
\begin{multline}\label{eq:eqdiff1}
\big(3z^3-z^4-3z^2+z\big)y''(z)+\big(5z^2\alpha-4z^3-2z^3\alpha+8z^2+1+\alpha-5z-4z\alpha\big)y'(z)\\
+\big(-1-2z^2-3z^2\alpha+2z-\alpha+4z\alpha-\alpha^2+2z\alpha^2-z^2\alpha^2\big)y(z)=0.
\end{multline}
This immediately translates into a linear recurrence satisfied by the sequence $(P_n(\alpha))$:
\begin{multline}\label{eq:rec1}
(n+3)(n+3+\alpha)P_{n+3}(\alpha)-(3n^2+4n\alpha+14n +\alpha^2+9\alpha+17)P_{n+2}(\alpha)
\\
+(3n+5+2\alpha)(n+2+\alpha) P_{n+1} (\alpha)-(n+2+\alpha)(n+1+\alpha)P_n(\alpha) 
=0
\end{multline}
with $P_0(\alpha)=\frac1\alpha$, $P_1(\alpha)=\frac{1+\alpha+\alpha^2}{\alpha(\alpha+1)}$ and $P_2(\alpha)=\frac{4+5\alpha+6\alpha^2+4\alpha^3+\alpha^4}{2\alpha(\alpha+1)(\alpha+2)}$.

\subsection{Proof of~\eqref{eq:subset2}} \label{preuvesubset2}

The proof is very similar to that of \cite{gvalues} so we skip the details. Let $(P_n,Q_n)$ be $E$-approximations of $\xi\in\C\etoile$. If $(P_n)$ has the first asymptotic behavior \eqref{eqtheoaebun} of Theorem~\ref{theoasypn}, then so does $(Q_n)$ with the same parameters $d$, $q$, $u$, $v$, and the sum is over the same non-empty finite set of $\te$. Therefore $\xi = \frac{g_\te \Gamma(-u_\te)}{g'_\te \Gamma(-u'_\te)} \in \Gamma(\Q)\cdot \Frac\GG$, using Eq. \eqref{eq123}.

Now if $(P_n)$ satisfies \eqref{eqtheoaebde} then so does $(Q_n)$ with the same parameters $q$, $u$, $v$, $\capa_1$, \ldots, $\capa_{d-1}$ (since we may assume that $d$ is the same), and the same set of $(\te_1,\ldots,\te_d)$ in the sum. If $v = \capa_1=\ldots=\capa_{d-1}=0$ and a term in the sum corresponds to $\te_1=\ldots=\te_{d-1}=0$, then $\xi = \frac{\omega_{0,\ldots,0,\te_d}}{\omega'_{0,\ldots,0,\te_d}}\in \frac{\E\cup(\Gamma(\Q)\cdot\GG)}{\E\cup (\Gamma(\Q)\cdot\GG)}$, else $\xi \in \Gamma(\Q ) \cdot \exp(\Qbar) \cdot \Frac\GG$ (using Eq. \eqref{eq123}).

\subsection{Extended $E$-approximations}\label{ssec:extended}

Let us consider the $E$-function
$$
E(z)=\sum_{n=1}^{\infty} \frac{z^n}{n!n}.
$$
We shall  prove that the sequence $(P_n)$ defined in the introduction by
$$
\frac{\log(1-z)}{1-z}-\frac{1}{1-z}E\left(-\frac z{1-z}\right) = \sum_{n=0}^{\infty} P_n z^n \in \mathbb Q[[z]]
$$
provides, together with $Q_n=1$, a sequence of $E$-approximations of Euler's constant in the extended sense of \eqref{eq:eapproxgen}. It is easy to see that
$$
P_n= \sum_{k=1}^n (-1)^{k-1} \binom{n}{k}\frac{1}{k!k}-\sum_{k=1}^n \frac1k = \sum_{k=1}^n (-1)^{k} \binom{n}{k}\frac{1}{k} \Big(1-\frac1{k!}\Big),
$$
where the second equality is a consequence of the identity $\sum_{k=1}^n \frac1k=\sum_{k=1}^n (-1)^{k-1} \binom{n}{k}\frac{1}{k}$.
We now observe that $E(z)$
 has the asymptotic expansion 
\begin{equation} \label{eqex2}
E(-z) \devasy - \gamma-\log(z) - e^{-z} \sum_{n=0}^{\infty} (-1)^n 
\frac{n!}{z^{n+1}}
\end{equation}
in a large sector bisected by any $\theta\in(-\pi,\pi)$ (see   \cite[Prop. 1]{Michigan}; this is also a special case of Theorem \ref{theoprecis}). 
 Therefore, for $z\to1$, $\vert z\vert <1$, 
$$
-\frac{1}{1-z} E\Big(-\frac{z}{1-z}\Big) + \frac{\log(1-z)}{1-z} =\frac{\gamma}{1-z} +   \mathcal{O}(1).
$$
As in \S \ref{ssec:example} in the case of $\Gamma(\alpha)$, a transfer principle readily shows that 
$$
\lim_{n\to +\infty} P_n = \gamma.
$$
Since $E(z)$ is holonomic, this is also the case of $\frac{\log(1-z)}{1-z}-\frac1{1-z} E\left(-\frac{z}{1-z}\right)$. The latter function satisfies the differential equation
\begin{multline}\label{eq:eqdiff2}
\big(3z^3- z^4-3z^2+z\big)y''(z)
+\big(1-5z+8z^2-4z^3\big)y'(z)+\big(-2z^2+2z-1\big)y(z)=0.
\end{multline}
This immediately translates into a linear recurrence satisfied by the sequence $(P_n)$:
\begin{equation}\label{eq:rec2}
(n+3)^2P_{n+3}- (3n^2+14n+17)P_{n+2}+(n+2)(3n+5)P_{n+1}-(n+1)(n+2)P_{n}=0
\end{equation}
with $P_0=0$, $P_1=0$, $P_2=\frac{1}{4}$. The differential equation~\eqref{eq:eqdiff2} and the recurrence relation~\eqref{eq:rec2} are the case $\alpha=0$ of  \eqref{eq:eqdiff1} and  \eqref{eq:rec1} respectively.

\bigskip

Let us now prove that any number with extended $E$-approximations is of the form \eqref{eq:subset3} stated in the introduction. 
Let $P(z)$ be given by \eqref{eq:eapproxgen}. If there is only one term in the sum, Theorems \ref{theo:eapprox} and \ref{theoasypn} hold and the proof extends immediately, except that $\E$ has to be replaced with $\E \cdot \log(\Qbaretoile)$ in \S \ref{subsec421} and \ref{subsec422}, and therefore in \eqref{eq:subset2} and \eqref{eq43}. Otherwise, we apply a variant of Lemma \ref{lemEsomme} to each $E$-function $E_\ell(z)$, obtaining exponential terms $e^{\alpha_{k,\ell} z }$: for each $k$ we write sufficiently many terms in the asymptotic expansion before the error term $o(1)$ (and not only the dominant one as in \S \ref{sec:asympPn}). Theorem \ref{thintrocce} asserts that all these terms are of the same form, but now the constants $\omega$ belong to ${\bf S}$. Combining these expressions yields
$$P(z) = \sum_{k=1}^K \omega_k e^{\alpha_k C(z)} U_k(z) (\log V_k(z))^{v_k}(1+o(1))$$
as $z $ tends to some point (possibly $\infty$) at which $C$ is infinite; here $U_k$, $V_k$ are algebraic functions, $v_k\in\N$, and $\omega_k\in{\bf S}$. However there is no reason why $\omega_k$ would belong to $\Gamma(\Q)\cdot \GG$ in general, since it may come from non-dominant terms in the expansions of $E_\ell(z)$, due to compensations. Upon replacing $\Gamma(\Q)\cdot \GG$ with ${\bf S}$ (and $\E$ with $\E \cdot \log(\Qbaretoile)$ as above), the proof of Theorems \ref{theo:eapprox} and \ref{theoasypn} extends immediately.

\bigskip

To conclude this section, we discuss another interesting example, which was also mentioned in the introduction. It corresponds to the more general 
notion of extended $E$-approximations where the coefficients of the linear form  \eqref{eq:eapproxgen} are in $\E$ and not just in $\Qbar$.
Let us consider the $E$-function $F(z^2)=\sum_{n=0}^{\infty} z^{2n}/n!^2$. It is solution 
of an $E$-operator $L$ of order $2$ with another solution of the form $G(z^2)+\log(z^2)F(z^2)$ where $G(z^2)=-2\sum_{n=0}^{\infty}\frac{1+\frac 12+\cdots +\frac 1n}{n!^2}z^{2n}$ is an $E$-function (in accordance with Andr\'e's theory). 
Then, 
$$
F(1-z)=\sum_{n=0}^{\infty} \frac{(1-z)^n}{n!^2}=\sum_{k=0}^{\infty} \frac{(-1)^k A_k}{k!} z^k 
$$
with 
\begin{equation*}\label{eq:Ak}
A_k=(-1)^k \sum_{n=0}^{\infty} \frac{1}{n!(n+k)!}.
\end{equation*}
It is a remarkable (and known) fact that the sequence $A_k$ satisfies the recurrence relation $A_{k+1}=kA_k+A_{k-1}$, $A_0=F(1), A_1=-F'(1)$. This can be readily checked. 
It follows that $A_{k}=V_kF(1)-U_kF'(1)$ where the sequences of integers $U_k, V_k$ are solutions of the same recurrence. 

Hence, the sequence 
$U_k/V_k$ is the sequence of convergents to $F(1)/F'(1)$ whose continued fraction is $[0;1,2,3,4,\ldots ]$. Moreover, we have
$$
\sum_{n=0}^{\infty} \frac{(-1)^k U_k}{k!}z^k = aF(1-z)+bG(1-z) + b \log(1-z) F(1-z)
$$
$$
\sum_{n=0}^{\infty} \frac{(-1)^k V_k}{k!}z^k = cF(1-z)+dG(1-z) + d \log(1-z) F(1-z)
$$
for some constants $a,b,c,d$, 
because both generating functions are solutions of 
an operator of order $2$ 
obtained from $L$ by changing $z$ to $\sqrt{1-z}$. The conditions $V_0=1, U_0=0, V_1=0,U_1=1$ and
$A_{k}=V_kF(1)-U_kF'(1)$ translate into a linear system in $a,b,c,d$ with solutions given by 
\begin{align*}
a&=\frac{g}{gf'-f^2-fg'}\in \E, \qquad b=-\frac{f}{gf'-f^2-fg'} \in \E,
\\
c&=-\frac{f+g'}{gf'-f^2-fg'} \in \E, \qquad d=\frac{f'}{gf'-f^2-fg'} \in\E,
\end{align*}
where $f=F(1), f'=F'(1), g=G(1), g'=G(1)$. We observe that $gf'-f^2-fg' \in\Qbar^*$ because it is twice the value at $z=1$ of the 
wronskian built on the linearly independent solutions $F(z^2)$ and $G(z^2)+\log(z^2)F(z^2)$. It follows that $U_k/V_k$ are extended $E$-approximations to the number $F(1)/F'(1)$ with ``coefficients'' in $\E$, but not in $\Qbar$ (because the number $f$ was proved to be transcendental by Siegel).

\def\refname{Bibliography}

S. Fischler, \'Equipe d'Arithm\'etique et de G\'eom\'etrie Alg\'ebrique, 
Universit\'e  Paris-Sud, B\^atiment 425,
91405 Orsay Cedex, France

T. Rivoal, Institut Fourier, CNRS et Universit\'e Grenoble 1, 
100 rue des maths, BP 74, 38402 St Martin d'H\`eres Cedex, France

\end{document}